\documentclass[review,hidelinks,onefignum,onetabnum]{siamart220329}

\usepackage{upgreek}
\usepackage{enumitem}
\usepackage{multirow}
\usepackage{wrapfig}

\usepackage{lipsum}
\usepackage{amsfonts}
\usepackage{graphicx}
\usepackage{epstopdf}
\usepackage{algorithmic}
\usepackage{csquotes}
\usepackage{float}
\usepackage[absolute,overlay]{textpos}
\allowdisplaybreaks
\ifpdf
  \DeclareGraphicsExtensions{.eps,.pdf,.png,.jpg}
\else
  \DeclareGraphicsExtensions{.eps}
\fi

\nolinenumbers
\newsiamremark{remark}{Remark}
\newsiamremark{hypothesis}{Hypothesis}
\crefname{hypothesis}{Hypothesis}{Hypotheses}
\newsiamthm{claim}{Claim}

\headers{Probable Event Constrained Optimization}{Qifeng Li}
\title{Probable Event Constrained Optimization\\ and A Data-embedded Solution Paradigm\thanks{Submitted to the editors DATE.
}}

\author{Qifeng Li\thanks{University of Central Florida, Orlando, FL 32816, USA 
  .}}

\usepackage{amsopn}

\ifpdf

\fi

\begin{document}

\maketitle

\begin{abstract}
This paper solves a new class of optimization problems under uncertainty, called Probable Event Constrained Optimization (PECO), which optimizes an objective function of decision variables and subjects to a set of Probable Event Constraints (PEC). This new type of constraint guarantees that optimal solutions are feasible for all uncertain events whose joint probabilities are greater than a user-defined threshold. The PEC can be used as an alternative to the conventional chance constraint, while the latter cannot guarantee the solution's feasibility to high-probability uncertain events. Given that the existing solution methods of optimization problems under uncertainty are not suitable for solving PECO problems, we develop a novel data-embedded solution paradigm that uses historical measurements/data of the uncertain parameters as input samples. This solution paradigm is conceptually simple and allows us to develop effective data-reduction schemes which reduce computational burden while preserving high accuracy. 

\end{abstract}



\section{Introduction}
\subsection{Background and motivation}
Uncertainties exist in the decision-making processes of many engineering systems. Generally, the procedure for solving decision-making problems under uncertainty consists mainly of three steps: 1) logic modeling, 2) deterministication\footnote{Inspired by the widely adopted term ``convexification" in the optimization field, which refers to the process of converting or approximating nonconvex problems into convex ones, we define a verb ``determinisfy" to refer to convert a logic model into its deterministic approximations, and the term ``deterministication" to refer to the process of determinisfying.}, and 3) solver solution. First, a logic model (or called uncertainty-perturbed model), such as stochastic (SO) \cite{birge2011introduction}, robust (RO) \cite{ben2009robust}, chance-constrained (CCO) \cite{li2008chance} optimizations, and their variants \cite{sahinidis2004optimization}, is chosen to mathematically formulate this decision-making problem under uncertainty. Second, the logic model is approximated by a deterministic program, i.e. linear, nonlinear, or integer program. Then, the resulting deterministic approximation is solved by off-the-shelf optimization algorithms \cite{nocedal2006numerical}/solvers \cite{anand2017comparative}. Generally, there are two basic requirements for the deterministic approximation, i.e., to be accurate and computationally tractable. 

In engineering fields, on one hand, engineers are generally interested in obtaining an optimal decision that works for all probable outcomes of uncertain events\footnote{See Definition \ref{def:pe}}. On the other hand, they are also interested in a solution paradigm that is applicable to nonlinear and even nonconvex problems since the mathematical models of many engineering systems, for example, power grids \cite{li2016convex}, gas/water supply networks \cite{li2018micro}, and transportation systems \cite{sharma2024enhanced}, are nonlinear and nonconvex. To meet the first engineering requirement, we propose a new logic model for mathematically formulating decision-making problems under uncertainty--the probabilistic event constrained optimization (PECO)--formulated as (\ref{PCCSP}). For the second, we propose a data-embedded solution paradigm, of which the resulting deterministic approximation is given in (\ref{DDA}). The proposed solution paradigm is less sensitive to the model's nonlinearity and nonconvexity than the Wald's \textit{minimax} paradigm (the latter one is widely adopted in existing work on optimization under uncertainty, for which a literature review is provided in Section 1.3).

\subsection{The proposed method}

Let $e$ denote an event in the event set $\mathcal{E}$ and $\mathbb{P}[\cdot]$ represent the probability of an event, this paper has the following definition on a ``probable event".
\begin{definition}[on a probable event] \label{def:pe}
    An event $e \in \mathcal{E}$ is called a probable event if $\mathbb{P}[e]\ge \alpha$, where $\alpha$ is a user-defined probability threshold.
\end{definition}
This research considers a random event of ``the outcome of the uncertain parameters\footnote{Called random variables in some references.} $\xi \in \Upxi$ equals to a specific value $y \in \Upxi$," where $\Upxi \subset (\mathbb{Z}^{u_1},\mathbb{R}^{u_2})$ ($u_1+u_2=u$), in the studied decision-making process. For example, in electric power systems, the operators need to deal with a random event that the power output of renewable energy resources are specific values. Then, considering $e:\xi=y$, i.e., the value of $\xi$ is $y$, as an event, we have the following specific formulation of PECO:
\begin{subequations} \label{PCCSP} 
\begin{align}
\text{\textbf{PECO}:}\quad\quad \min_{x\in \mathbb{R}^n} \quad & f(x)  \label{obj_GOU} \\
 \mathrm{s.t.} \quad  & g(x,y) \le 0,\;\forall y \in \Upxi_\alpha = \{ y \in \Upxi \,|\,\mathbb{P}[\xi=y] \ge \alpha\} \label{PCC} 
\end{align} 
\end{subequations} 
where $f:\mathbb{R}^n \rightarrow \mathbb{R}$, and $g:\mathbb{R}^n \times \Upxi \rightarrow \mathbb{R}^m$ are continuous and differentiable functions. Without loss of generality, only inequality constraints are considered in (\ref{PCC}) since there are explicit and implicit methods of equivalently reformulating equations as inequations. The vector of uncertain parameters $\xi$ follows some probability distributions that may be \textit{unknown} and $y$ represents a future outcome of $\xi$. 

Note that an outcome of the uncertain parameters is also called a realization/ observation in the literature. This paper considers two types of outcomes for the uncertain parameter $\xi$, i.e., the ``scenario" and ``data point," which have different meanings from each other. In the rest of the paper, the proposed research is narrated based on the terminologies defined below.
\begin{definition}[on key terminologies] \label{def:srd}
\begin{itemize}[leftmargin=0.36cm]
    \item \textbf{Data point} ($\xi_{\rm d}$): the value of $\xi$ that is measured in history, which can be considered an independent and identically distributed (i.i.d.) sample. 
    \item \textbf{Scenario} ($\xi_{\rm s}$): a possible outcome of $\xi$ that a probability is assigned to event $\xi=\xi_{\rm s}$, which can represent a set of data points.
    \item \textbf{Probable scenario}: a scenario of $\xi$ that satisfies $\mathbb{P}[\xi=\xi_{\rm s}] \ge \alpha$,  where $\alpha$ is a user-defined probability threshold. 
    \item \textbf{Probable data point}: a data point $\xi_{\rm d}$ of $\xi$ that satisfies $\xi_{\rm d}=\xi_{\rm s}$ is called a probable data point, where $\xi_{\rm s}$ is a probable scenario.
     \item \textbf{Scenario set} $\mathcal{S}=\{\xi^{(k)}_{\rm s},\; k=1,\ldots,S\}$: is a set of scenarios where each of its elements is unique, i.e., $\xi^{(i)}_{\rm s} \neq \xi^{(j)}_{\rm s}$ if $i \neq j$ ($1 \le i,j\le S$). 
    \item \textbf{Data set} $\mathcal{D}=\{\xi^{(k)}_{\rm d},\; k=1,\ldots,D\}$: a \textbf{finite multiset} of data points of $\xi$, where there may be multiple instances for each of its elements, that is, $\xi^{(i)}_{\rm d} = \xi^{(j)}_{\rm d}$ is possible even if $i \neq j$ ($1 \le i,j\le S$). Note that $D$ is generally large. 
    \item $\Omega$ denotes the\textbf{ mathematical representation} of the uncertainty set\footnote{Called sample space in some references.} of $\xi$, namely $\Omega$ uses mathematical constraints to describe $\Upxi$. 
    \item $\mathcal{S}^\forall$ denotes the \textbf{particle\footnote{A particle may mean a scenario, sample, or data point in different contexts.} representation} of $\Upxi$, which is a scenario set that contains all the possible scenarios of $\xi$. $\mathcal{S}^\forall$ is generally \textbf{infinite} for continuous or mixed-integer $\xi$.
\end{itemize}    
\end{definition}
Moreover, this paper uses $\Omega$ ($\Omega_\alpha$) and $\mathcal{S}^\forall$ ($\mathcal{S}^\forall_\alpha$) to denote the mathematical and particle representations of $\Upxi$ ($\Upxi_\alpha$), respectively. In other words, given in (\ref{PCC}) that $\Upxi_\alpha = \{ y \in \Upxi \,|\,\mathbb{P}[\xi=y] \ge \alpha\}$, we have
\begin{subequations}
    \begin{align}
    &\Omega_\alpha=\{ y \in \Upxi \,|\,P(y) \ge \alpha\}\,(\text{if probability distribution function}\, P\, \text{is given}), \label{omegaalpha} \\
    &\mathcal{S}^\forall_\alpha=\{\forall \xi^{(k)}_{\rm s} \in \mathcal{S}^\forall \,|\, \mathbb{P}[\xi=\xi^{(k)}_{\rm s}] \ge \alpha \}.
    \end{align}
\end{subequations}

 Assuming that $P(\xi)$ is not perfectly known, and instead, a set of historical data $\mathcal{D}$ is available, this paper proposes a data-embedded deterministic approximation (DeDA) of PECO (\ref{PCCSP}) which is in the form of
\begin{subequations} \label{DDA}
\begin{align}
\text{\textbf{DeDA}:}\quad\quad \min_{x\in \mathbb{R}^n} \quad & f(x)  \nonumber \\
 \mathrm{s.t.} \quad  &  g(x,\xi_{\rm d}^{(k)}) \le 0,\;   ( \xi_{\rm d}^{(k)}  \in \Box \subset\mathcal{D}_\alpha) 
\end{align} 
\end{subequations}
where $\mathcal{D}_\alpha$ is the subset of $\mathcal{D}$ that contains all probable data points, and $\Box$ represents a subset of $\mathcal{D}_\alpha$ that is selected according to specific criteria which will be discussed in the rest of the paper.

Although the logic model (\ref{PCCSP}) looks similar to some existing models like the RO approximation of CCO (RO-CCO) \cite{ben2009robust} and distributionally robust optimization (DRO) \cite{rahimian2019distributionally}, we will show later that the PECO (\ref{PCCSP}) is different from RO-CCO and DRO in terms of both logical meanings and some key mathematical properties. It's also worth noting that the term ``data-embedded" means the direct substitution of uncertain variables with data points in the deterministic approximation, which is different from the ``data-driven" methods that we will review in the next subsection. 

\subsection{Related work}
Although the PECO was proposed for engineering needs, it is somewhat related to existing work on optimization under uncertainty. Given the presence of the sign ``$\forall$" and the data-based nature of the proposed method, we considered that, among the existing work, those on RO-CCO (\ref{ROCCO}) and the DRO (\ref{DRO}) are the most relevant to the proposed work. While the classic RO (\ref{RO}) may be overly-conservative, the classic CCO (\ref{CCO}) can be considered a less-conservative modeling paradigm for optimization under uncertainty. For linear cases, it's often approximated by the RO-CCO, of which a typical formulation is given in (\ref{ROCCO}), then solved via the computationally effective Wald's \textit{minimax} paradigm (\ref{DROCCO}) \cite{sniedovich2008wald}.
\begin{align}
    \text{\textbf{Classic RO}:}\quad\min_{x\in \mathbb{R}^n} \; & f(x) && \mathrm{s.t.}\; g(x,\xi) \le 0, \; \forall \xi \in \Upxi \label{RO}\\
      \text{\textbf{Classic CCO}:}\quad\min_{x\in \mathbb{R}^n} \; & f(x) && \mathrm{s.t.}\; \mathbb{P}[g(x,\xi) \le 0] \ge 1-\beta \label{CCO} \\
    \text{\textbf{RO-CCO}:}\quad\min_{x\in \mathbb{R}^n} \; & f(x) && \mathrm{s.t.}\; g(x,\xi) \le 0, \; \forall \xi \in \mathcal{U}_\beta \label{ROCCO}\\
     \textbf{minimax}_{\text{RO-CCO}}:\quad\min_{x\in \mathbb{R}^n} \; & f(x)&& \mathrm{s.t.}\; \max_{\xi \in \mathcal{U}_\beta} \{g(x,\xi)\} \le 0,  \label{DROCCO}
\end{align}
That's to say, if one can construct an uncertainty set $\mathcal{U}_\beta$ that satisfies
\begin{equation} \label{uset}
   \mathcal{U}_\beta = \arg_{\mathcal{U}} \left\{ \int\cdots\int_{\mathcal{U}}P(\xi)d\xi_u\cdots d\xi_1\ge 1-\beta \right\},
\end{equation}
where $\arg$ denotes the argument of a function and recall that $P(\cdot)$ is the probability distribution function of $\xi$, one can guarantee that a solution of RO-CCO (\ref{ROCCO}) is feasible to CCO (\ref{CCO}). However, there are some limitations associated with the modeling and solution paradigm of RO-CCO. First, in reality, there may exist multiple $\mathcal{U}_\beta$'s that satisfy condition (\ref{uset}), which means that (\ref{ROCCO}) is not an equivalent to (\ref{CCO}) even a proper set $\mathcal{U}_\beta$ is found. 
Second, the \textit{minimax} model (\ref{DROCCO}), which was originally designed for linear cases, may be computationally intractable or even unable to yield a meaningful solution for nonlinear, nonconvex problems. As a result, the existing research basically solved linear RO-CCO problems with a focus on constructing a linear or convex uncertainty set $\mathcal{U}_\beta$, such that the resulting \textit{minimax}-based deterministic approximations possess computational efficiency \cite{ben2011chance,margellos2014road,li2015optimal,yuan2017robust,han2016robust,bertsimas2004price,chen2010cvar}.

More often than not, decision makers know only partial information on $P$, that is, we know only that $P$ belongs to a given family $\mathcal{P}$ of distributions. DRO (\ref{DRO}) is a good fit for modeling the corresponding optimization problems under uncertainty, of which a typical formulation is given in (\ref{DRO}) where the uncertainty set $\mathcal{P}$ is also called the ambiguity set. DRO can be considered a combination of the classic RO and the classic SO. First, it can be considered an SO problem of which the solutions are robust to an ambiguity set, i.e., a set of probability distributions. Second, it can be considered an RO problem where the uncertainty is a probability distribution and the corresponding uncertainty set is the ambiguity set.
\begin{align}
        \text{\textbf{DRO}:}\quad\quad\min_{x\in \mathbb{R}^n} \; & f(x) && \mathrm{s.t.}\; \mathbb{E}[g(x,\xi)] \le 0, \; \forall P \in \mathcal{P} \label{DRO} \\
        \textbf{minimax}_{\text{DRO}}:\quad\min_{x\in \mathbb{R}^n} \; & f(x) && \mathrm{s.t.}\; \max_{P \in \mathcal{P}} \{g(x,\xi)\} \le 0.  \label{DDRO}
\end{align}

Being similar to RO-CCO, the most common solution process of DRO also involves approximating the logic model (\ref{DRO}) with a \textit{minimax}-based deterministic model (\ref{DDRO}). The main difference lies in that the decision variables of the lower-level optimization problem of (\ref{DDRO}) are the probability distribution functions $P$. Then, the existing research mainly focused on constructing ambiguity sets with a goal of making (\ref{DDRO}) computationally tractable, e.g., the moment-based \cite{calafiore2006distributionally,delage2010distributionally,hanasusanto2015distributionally}, statistical distance-based \cite{hu2013kullback,mohajerin2018data}, and likelihood-based \cite{wang2016likelihood} ambiguity sets. 
Although the DRO paradigm allows decision makers to incorporate partial distribution information learned from uncertainty data into the optimization, it also subjects to the limitations of the \textit{minimax}-based solution method.

\subsection{Contributions}
The research presented in this paper is different from the existing research on the following three aspects. First, while the RO-CCO is an approximations of the classic CCO, PECO is an alternative to CCO due to its different logical meaning. The PEC (\ref{PCC}) logically means that ``a solution of problem (\ref{PCCSP}) should satisfy $g(x,\xi) \le 0$ for all probable outcomes of the uncertain parameters $\xi \in \Upxi$" (denoted as logical meaning (LM)$\#$1). In contrast, the logical meaning of the conventional chance constraint is that ``the probability of that a solution violates constraint $g(x,\xi) \le 0$ is not bigger than $\beta$" (denoted as LM$\#$2). It's worth noting that PEC also guarantees LM2 under certain condition while the chance constraint does not guarantee LM1. That means PECO is a ``safe" alternative to the CCO. If the probability distribution function $P(\cdot)$ of $\xi$ is known a priori, the mathematical difference between PECO and RO-CCO mainly resides in the difference between $\Omega_\alpha$ (\ref{omegaalpha}) and $\mathcal{U}_\beta$ (\ref{uset}).


Second, the proposed deterministic approximation of PECO, i.e., the DeDA (\ref{DDA}) is not \textit{minimax}-based. In the existing research (as we reviewed in the previous subsection), the solution process of the RO-CCO and DRO paradigms is generally based on the \textit{minimax} model which is a bi-level optimization problem. Solving such a bi-level optimization model for nonlinear, nonconvex problems, on one hand, is computationally difficult. On the other hand, the obtained solution may not be meaningful to the original RO-CCO or DRO problem since there is no theoretical guarantee for identifying the ``worst case" due to the presence of multiple local maxima for the lower-level problem. In contrast, the proposed DeDA (\ref{DDA}) does not subject to these limitations.

f which the solution process rely heavily on a structurally simple subproblem, i.e., linear or at least convex problems. As a result, the existing research mainly focused on constructing a linear or convex $\mathcal{U}_\beta$ ($\mathcal{P}$) for RO-CCOs (DROs). In contrast, the proposed solution paradigm does not involve the construction of $\mathcal{U}_\beta$ or $\mathcal{P}$ due to its data-based nature.

Third, the proposed solution paradigm (\ref{DDA}) is ``data-embedded" rather than ``data-driven." As being pointed out in Subsection 1.2, the existing research on RO-CCO and DRO mainly focused on constructing $\mathcal{U}_\beta$ or $\mathcal{P}$, of which the existing methods are basically: 1) distribution-driven, assuming that the probability distributions of uncertain parameters are known a priori \cite{ben2011chance,ben2009robust,bertsimas2004price}; or 2) data-driven, when the distributions are not perfectly known and a set of historical data is available instead \cite{delage2010distributionally,ben2013robust,bertsimas2018data}. For examples, a systematic method, that uses hypothesis test to construct uncertainty sets from data, was presented in \cite{bertsimas2018data}, and a data-driven ambiguity set was proposed in \cite{delage2010distributionally} based on the distribution’s support information as well as the confidence regions for the mean and second-moment matrix. While the existing data-driven methods mainly use historical data to construct $\mathcal{U}_\beta$ or $\mathcal{P}$ (i.e., infer a mathematical formulation of $\mathcal{U}_\beta$ or $\mathcal{P}$), the proposed data-embedded method directly input data points into DeDA (\ref{DDA}) instead of using these data to construct uncertain models, such as $\mathcal{U}_\beta$ or $\mathcal{P}$. 

Given the above differences, the contributions of this research can be summarized as follows: 1) the proposed PECO provides engineers a new, practical option of solving decision-making problems under uncertainty due to its new logical meaning that has a high practical value in engineering; 2) it's less sensitive to nonlinearity and nonconvexity since the proposed deterministic approximation of PECO is not \textit{minimax}-based; 3) the proposed data-embedded solution paradigm does not involve uncertainty model construction which is a bottleneck in the existing research; and 4) algorithms of strategic data selection are developed for improving the computational efficiency of the proposed solution paradigm by effectively eliminating inactive data points (see Section 4).

\section{The Probable Event Constraint (PEC)}

This section aims at revealing more properties of PEC (\ref{PCC}) and, more importantly, showing that the existing approaches of solving optimization problems under uncertainty may not be applicable for solving PECO.

 \subsection{Illustrative examples of Definition \ref{def:srd}}
\begin{wrapfigure}{r}{0.3\textwidth}
\includegraphics[width=3.8cm]{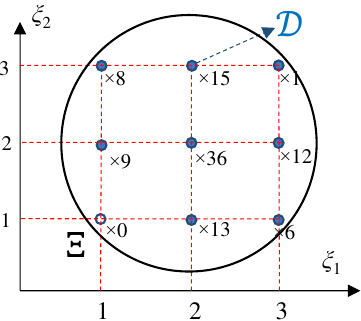}
\vspace{-6pt}
\caption{An illustrative example of Definition \ref{def:srd}.}
\label{fig:realization_scenario}
\vspace{-12pt}
\end{wrapfigure}
Recall that this paper considers two types of outcomes for the uncertain parameter $\xi$, i.e., the scenarios and data points,
Figure \ref{fig:realization_scenario} provides an illustrative example for the terminologies defined in Definition \ref{def:srd}, where the uncertain vector $\xi =[\xi_1,\xi_2]^{\rm T}$ and $\Omega=\{\xi  \in \mathbb{Z}^2 | \xi_1^2+\xi_2^2-4\xi_1-4\xi_2\le -5.9\}$ is the mathematical representation of the uncertainty set $\Upxi$. This example considers 100 data points that are measured in history, namely the data set $\mathcal{D}=\{\xi^{(k)}_{\rm d},\; k=1,\ldots,100\}$. If all the scenarios are ordered as $\xi^{(1)}_{\rm s}=(1,1)$, $\xi^{(2)}_{\rm s}=(1,2)$, $\xi^{(3)}_{\rm s}=(1,3)$, $\xi^{(4)}_{\rm s}=(2,1)$, $\xi^{(5)}_{\rm s}=(2,2)$, $\xi^{(6)}_{\rm s}=(2,3)$, $\xi^{(7)}_{\rm s}=(3,1)$, $\xi^{(8)}_{\rm s}=(3,2)$, and $\xi^{(9)}_{\rm s}=(3,3)$, we have $\mathcal{S}^\forall=\{\xi^{(k)}_{\rm s},\; k=1,\ldots,9\}$ which is a finite set since $\xi_1$ and $\xi_2$ are two finite integer parameters. From Figure \ref{fig:realization_scenario}, we know that the scenario $\xi^{(2)}_{\rm s}$ occurred 9 times in the data set $\mathcal{D}$, and hence we can let $\xi^{(j)}_{\rm d}=\xi^{(2)}_{\rm s}$ ($j=1,\ldots,9$). If we consider the joint probability of $\xi^{(2)}_{\rm s}$ as $\mathbb{P}[\xi =\xi^{(2)}_{\rm s}]=9/100=0.09$, the probabilities of other scenarios can be obtained in the same way. If we set $\alpha=0.1$, the set of probable scenarios $\mathcal{S}^\forall_\alpha=\{\xi^{(k)}_{\rm s},\; k=4, 5, 6, 8\}$. For this particular example, $\Omega_\alpha$ (i.e., the mathematical expression of $\Upxi_\alpha$) may not exist. $\mathcal{D}_\alpha$ contains all data points that equal to $\xi^{(k)}_{\rm s}$ ($k=4, 5, 6, 8$) and $|\mathcal{D}_\alpha|=85$. It's not hard to know that, for this example, $\mathcal{S}^\forall_\alpha=\mathfrak{U}[\mathcal{D}_\alpha]$ where $\mathfrak{U}[\cdot]$ denotes the underlying set\footnote{A underlying set is the set of distinct elemenets of a multiset} of a data set. 

\begin{wrapfigure}{r}{0.3\textwidth}
\includegraphics[width=3.8cm]{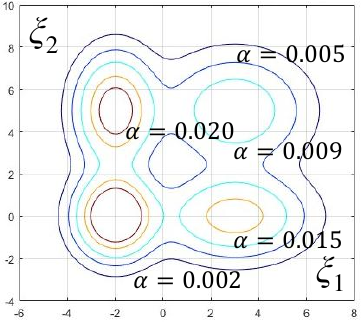}
\vspace{-6pt}
\caption{An illustrative example of $\Upxi_\alpha$.}
\label{fig:contour}
\vspace{-12pt}
\end{wrapfigure}
 An example of $\Upxi_\alpha$ is provided in Figure \ref{fig:contour}, where $\xi_1$ and $\xi_2$ follow bimodal distributions. To be specific, $\mathbb{P}[\xi_1=y_1] \sim [0.5\mathcal{N}(-2,1)+0.5\mathcal{N}(3,1.8)]$ and $\mathbb{P}[\xi_2=y_2] \sim [0.5\mathcal{N}(0,1.2)+0.5\mathcal{N}(5,1.6)]$, where $\mathcal{N}(\mu,\sigma^2)$ is a normal distribution with $\mu$ as the mean and $\sigma^2$ as the variance. Although, for this specific example, it's not hard to obtain the mathematical expression of $\Upxi_\alpha$, i.e., $\Omega_\alpha$ as given in (\ref{upxialpha}), it is highly nonconvex and even becomes discontinuous as $\alpha$ varies. From this example, we observe that, in addition to the nonconvexity of system constraints, $\Upxi_\alpha$ can be highly nonconvex. What is worse, obtaining $\Omega_\alpha$ is not necessarily easy for many real-world cases, regardless of whether it is non-convex or convex. Hence, this paper proposes a non-\textit{minimax}-based, data-embedded deterministic approximation (\ref{DDA}) for PECO which does not rely on knowing $\Omega_\alpha$ (details are presented in Section 3.).
 \begin{equation}\label{upxialpha}
 \begin{aligned}
    \Omega_\alpha=&\{(\xi_1,\xi_2)\in \Upxi\,|\, (\frac{0.5}{\sqrt{2\pi}}e^{-(\xi_1+2)^2} \\
     &+\frac{0.5}{\sqrt{3.6\pi}}e^{-\frac{(\xi_1-3)^2}{1.8}})(\frac{0.5}{\sqrt{2.4\pi}}e^{-\frac{\xi_2^2}{1.2}}+ \frac{0.5}{\sqrt{3.2\pi}}e^{-\frac{(\xi_2-5)^2}{1.6}})\ge \alpha \}.
\end{aligned}
 \end{equation}
 

\subsection{Relations between PEC and the classic chance constraint} 
First, while the PECO (\ref{PCCSP}) can be considered an alternative to the CCO (\ref{CCO}), it's worth noting that the PEC (\ref{PCC}) has a different logic meaning from the conventional chance constraint. Let $x^*_{\rm P}$ denote the optimal solution of PECO (\ref{PCCSP}) and $\check{\xi}$ be an arbitrary outcome of $\xi \in \Upxi$, the logical meaning of the PEC (\ref{PCC}) is: ``$g(x^*_{\rm P},\check{\xi}) \le 0$\textit{ holds if $\check{\xi}$ is a probable outcome of $\xi$.}"
Further let $x^*_{\rm C}$ denote the optimal solution of the classic CCO (\ref{CCO}), we recall the following logical meaning of chance constraint for comparison purpose:
``\textit{The probability of} $g(x^*_{\rm C},\check{\xi}) \le 0$ \textit{should not be less than} $1-\beta$ \textit{(or, the probability of} $g(x^*_{\rm C},\check{\xi}) > 0$ \textit{should not be bigger than $\beta$)}."

Since the PECO (\ref{PCCSP}) and CCO (\ref{CCO}) share the same objective function, the difference between them lies in the constraints. 
If the mathematical expression of $P(\xi)$ is known, the feasible space of $x$ specified by the PEC (\ref{PCC}) is given as
\begin{equation} \label{feasisetpccsp}
    \mathcal{X}_{\rm P}=\{ x \in \mathbb{R}^n \mid g(x,y) \le 0,\;\forall y \in \Omega_\alpha\},
\end{equation}
where $\Omega_\alpha$ is in (\ref{omegaalpha}), and we have the following proposition.
\begin{proposition} \label{pro:feasisetpccsp}
    $\mathcal{X}_{\rm P}(\alpha_1) \supseteq \mathcal{X}_{\rm P}(\alpha_2)$ if $\alpha_1 \ge \alpha_2$.
\end{proposition}
The proof of this proposition is given in Appendix A.1. Recall the logical meaning of chance constraints, if $P(\xi)$ is known, a deterministic formulation of the feasible space of $x$ that is specified by the chance constraint in (\ref{CCO}) is given as 
\begin{equation} \label{feasisetcccsp1}
    \mathcal{X}_{\rm C}=\left\{ x \in \mathbb{R}^n \middle| 
\int\cdots\int_{ \mathcal{M}(x)}P(\xi)d\xi_u\cdots d\xi_1\ge 1-\beta \right\},
\end{equation}
where $\mathcal{M}(x):=\{\xi \in \Omega|g(x,\xi) \le 0, x \in \mathbb{R}^n\}$ is the set of $\xi$ that $g(x,\xi) \le 0$ holds for a specific $x$. We have the following propositions. 
\begin{proposition}[on the relations between $\mathcal{X}_{\rm P}$ and $\mathcal{X}_{\rm C}$] \label{pro:CD}
     If
    \begin{equation} \label{CD}
  \alpha = \arg_v \left\{ \int\cdots\int_{\{\xi \in \Omega|P(\xi)\ge v\}}P(\xi)d\xi_u\cdots d\xi_1=1-\beta \right\} ,
\end{equation}
where $\arg$ means the argument of a function, we have the following relations:
\begin{enumerate}[leftmargin=*]
    \item  $\mathcal{X}_{\rm P} \subseteq \mathcal{X}_{\rm C}$;
    \item  $\mathcal{X}_{\rm P} = \mathcal{X}_{\rm C}$ if, when $\xi^{(a)}$ and $\xi^{(b)}$ are arbitrarily realizations in $\mathcal{M}(x)$ and $\Omega \setminus \mathcal{M}(x)$ respectively, $P(\xi^{(a)}) \ge P(\xi^{(b)})$.
\end{enumerate}
\end{proposition}
The proof of Proposition \ref{pro:CD} is provided in Appendix A.2. Although the classic chance constraint restricts that its feasible solutions can ensure a satisfactory probability of constraint violation, it cannot guarantee that its optimal solution is feasible to a probable realization of $\xi$, which is not desirable in engineering. An advantage of PEC is that it can guarantee both. 

\subsection{Applying the scenario method to solve PECO}
Given that the scenario method \cite{campi2008exact,calafiore2010random,garatti2022risk,campi2008exact,calafiore2010random} is a distribution-free approach for solving CCO problems, this subsection evaluates its applicability for solving PECO. If the scenario method is directly applied, the PECO is approximated by 
\begin{subequations} \label{DA2}
\begin{align}
 \min_{x} \quad & f(x)    \label{obj_DA2} \\
 \mathrm{s.t.} \quad  &g(x,\xi^{(k)}_{\rm iid}) \le 0,\;   (k=1,\ldots,N)\label{constr_DA2} 
\end{align} 
\end{subequations}
where $\xi^{(k)}_{\rm iid}$ is an i.i.d. sample. Recall that $\mathcal{X}_{\rm P}$ is the feasible set of PECO and further let $\mathcal{X}_{\rm S}=\{ x \in \mathbb{R}^n \mid (\text{\ref{constr_DA2}) holds}\}$ denote the feasible set of the scenario-based deterministic approximation (\ref{DA2}) of PECO, we have the following proposition:
\begin{proposition} \label{RDA}
    $\mathcal{X}_{\rm P} \subset \mathcal{X}_{\rm S}$ ($\mathcal{X}_{\rm P} \neq \mathcal{X}_{\rm S}$) is a high-probability event if $N \le \frac{1}{\alpha}$ ($N > \frac{1}{\alpha}$).
\end{proposition}
Readers can find the proof of this proposition in Appendix A.3. When, in CCO (\ref{CCO}), $f(\cdot)$ is linear and $g(\cdot)$ is convex on $x$, let $x^*_{\rm C}$ and $x^*_{\rm S}$ denote the optimal solutions of CCO (\ref{CCO}) and its scenario-based deterministic approximation (\ref{DA2}) respectively, and $\varepsilon=\mathbb{P}[x^*_{\rm S}=x^*_{\rm C}]$, Theorem 1 in \cite{campi2008exact} asserts that $\varepsilon = \epsilon$ if
\vspace{-6pt}
\begin{equation} \label{N-RCP}
\vspace{-6pt}
   N \ge  \frac{e(n-\ln{\epsilon})}{\beta(e-1)}, 
\end{equation}
where $e$ is Euler's number. Nevertheless, the $N$ that satisfies (\ref{N-RCP}) is significantly bigger than $1/\alpha$ under condition (\ref{CD}). According to Proposition \ref{RDA}, it cannot guarantee that $\mathcal{X}_{\rm S}$ is an accurate approximation of $\mathcal{X}_{\rm P}$ with the $N$ determined by (\ref{N-RCP}).

In summary, the existing solution methods for optimization problems under uncertainty have limited applicability for solving PECO (i.e., the scenario approach is not applicable, while the minimax method is only applicable to some simple convex cases). The purpose of next section is developing a novel solution method that is: 1) distribution-free (it's based on historical data rather than pre-known probability distribution functions), 2) general (it does not rely on linearity or convexity assumptions on system constraints), and 3) computationally effective (it facilitates the development of effective data reduction methods).

\section{The Proposed Solution Paradigm for PECO: Data-embedded Deterministication} 
 This paper assumes that the probability distribution function $P(\xi)$ is not perfectly known, and instead, a set of historical data $\mathcal{D}$ as defined in Definition \ref{def:srd} is available. As mention in Subsection 1.2, we propose a data-embedded deterministic approximation, which is in the form of (\ref{DDA}), as the deterministication step of solving PECO problems. This subsection aims to reveal more details and properties of (\ref{DDA}) and the concept of data-embedded deterministication.

 \subsection{DeDA($\mathcal{D}_{\alpha}^z$): A data-embedded deterministic approximation for PECO} Noting that this paper uses DeDA ($\Box$) to denote the data-embedded optimization model (\ref{DDA}) with $\Box$ as the embedded data set, we first have the following proposition. 
\begin{proposition}\label{pro:underlying}
    Given two data sets $\mathcal{D}_1$ and $\mathcal{D}_2$, where $\mathcal{D}_2$ is a multiset. DeDA ($\mathcal{D}_1$) is equivalent to DeDA ($\mathcal{D}_2$) if $\mathcal{D}_1=\mathfrak{U}[\mathcal{D}_2]$, recalling that $\mathfrak{U}[\cdot]$ denotes the underlying set of a data set.
\end{proposition}
The proof of the above proposition can be found in Appendix B.1. Recalling that $\mathcal{D}_\alpha$ is the subset of $\mathcal{D}$ that contains all probable data points, Subsection 3.3 describes an approach for obtaining $\mathcal{D}_\alpha$ from $\mathcal{D}$. 

Second, let $\mathcal{D}_\alpha^z \subset \mathcal{D}_\alpha$ be a set of $z$ data points that are \textit{randomly} selected from $\mathcal{D}_\alpha$, we consider DeDA($\mathcal{D}_{\alpha}^z$) as the first deterministic approximation of PECO (\ref{PCCSP}). Next, denoting
\begin{equation} \label{FS}
    \mathcal{X}_{\rm D}(\Box)=\{ x \in \mathbb{R}^n \mid g(x,\xi_{\rm d}^{(k)}) \le 0\;   (\forall \xi_{\rm d}^{(k)}  \in \Box)\}
\end{equation}
as the feasible space of DeDA (\ref{DDA}) with $\Box$ as the embedded data set, we discuss the relations between the feasible spaces of DeDA($\mathcal{D}^z_{\alpha}$) and PECO in what follows. Recalling $\mathcal{S}_\alpha^{\forall}$ which was defined in Definition \ref{def:srd}, we have the following definitions.
\begin{definition}[on a deterministic equivalent of PECO]\label{def:equivalent}
    DeDA($\mathcal{S}_\alpha^{\forall}$) is considered a deterministic equivalent to PECO, i.e., $\mathcal{X}_{\rm D}(\mathcal{S}_\alpha^{\forall})=\mathcal{X}_{\rm P}$.
\end{definition}
\begin{definition}[on the boundary-forming data points of the feasible space of DeDA] \label{def:adpfs}
    If a data set $\mathcal{B}_{\Box}^{\rm FS}=\{\xi_{\rm d}^{(k)},\,k=1,\ldots,B_{\Box}^{\rm FS} \}$ is the SMALLEST subset of  $\Box$ that satisfies:
    \begin{equation} \label{condition_adpfs}
        \mathcal{X}_{\rm D}(\mathcal{B}_{\Box}^{\rm FS}) = \{ x \in \mathbb{R}^n \mid g(x,\xi_{\rm d}^{(k)}) \le 0\;   (\forall \xi_{\rm d}^{(k)} \in \mathcal{B}_{\Box}^{\rm FS}) \}=\mathcal{X}_{\rm D}(\Box),
    \end{equation}
    the elements in $\mathcal{B}_{\Box}^{\rm FS}$ are the boundary-forming data points of feasible space $\mathcal{X}_{\rm D}(\Box)$.
\end{definition}
  Note that, following expression (\ref{FS}), $\mathcal{X}_{\rm D}(\mathcal{S}_\alpha^{\forall})$ denotes the feasible space of DeDA($\mathcal{S}_\alpha^{\forall}$), where $\mathcal{S}_\alpha^{\forall}$ is not a multiset and $|\mathcal{S}_\alpha^{\forall}|=\infty$ when $\xi$ is continuous or mixed-integer. Following Definition \ref{def:adpfs}, the boundary-forming data set of $\mathcal{X}_{\rm D}(\mathcal{S}_\alpha^{\forall})$ is denoted as $\mathcal{B}_{\forall,\alpha}^{\rm FS}$. For the general case of $\xi$ (i.e., $\xi$ can be integer, continuous, or mixed-integer), we have the following proposition.

\begin{proposition}[on the relations between the feasible spaces of DeDA($\mathcal{D}_\alpha^z$) and PECO] \label{pro:adpfs}
    $\mathcal{X}_{\rm D}(\mathcal{D}_\alpha^z) \supseteq \mathcal{X}_{\rm P}$ and, if and only if $\mathcal{B}_{\forall,\alpha}^{\rm FS} \subseteq \mathcal{D}_\alpha^z$, $\mathcal{X}_{\rm D}(\mathcal{D}_\alpha^z) = \mathcal{X}_{\rm P}$.
\end{proposition}

\begin{wrapfigure}{r}{0.456\textwidth}
\includegraphics[width=5.8cm]{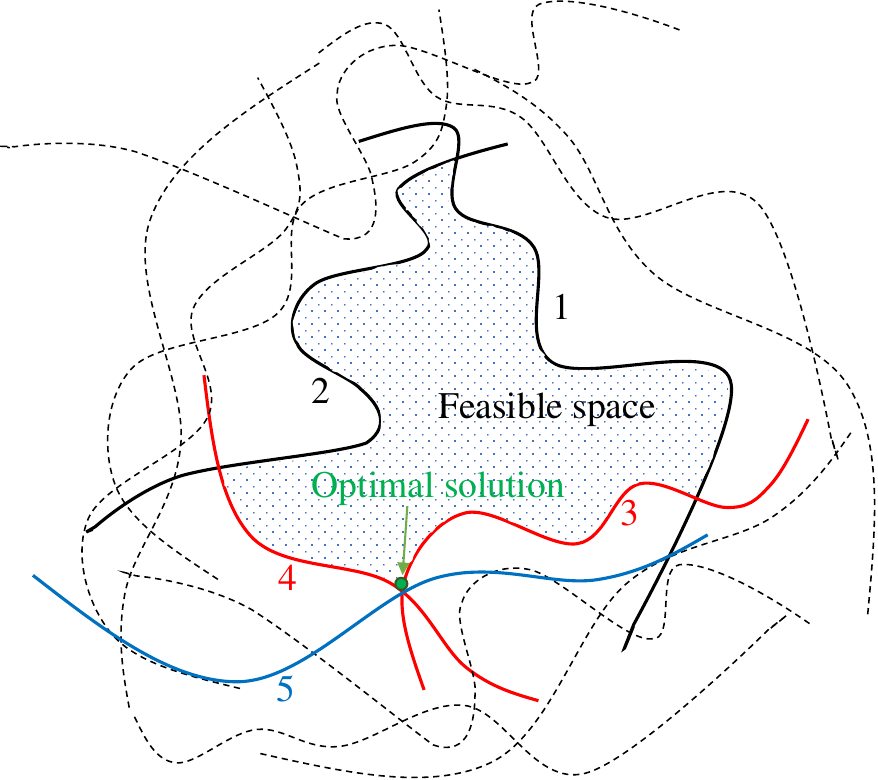}
\vspace{-6pt}
\caption{A pictorial interpretation of boundary-forming, active, and inactive data points in a 2-D $x$-space, where constraints 1-4 are the boundary-forming constraints of the feasible space, 3 and 4 are the boundary-forming constraints of the optimal solution, 3-5 are active constraints, and the rest are inactive constraints. A boundary-forming/active data point is a data point that contributes at least one boundary-forming/active constraint.}
\label{fig:activedata}
\vspace{-12pt}
\end{wrapfigure}
Readers can find the proof of this proposition in Appendix B.2. A pictorial explanation of definition \ref{def:adpfs} and proposition \ref{pro:adpfs} is provided in Figure \ref{fig:activedata}, which implies that the feasible space of a DeDA for a PECO problem is determined by a limited number of boundary-forming data points. The rest data points are inactive which can be removed without impacting the feasible space. Proposition \ref{pro:adpfs} indicates that, if one can guarantee that the finite data set $\mathcal{D}_\alpha^z$ contains all boundary-forming data points of the deterministic equivalent of PECO, i.e., the DeDA($\mathcal{S}_\alpha^{\forall}$), the finite deterministic optimization problem DeDA($\mathcal{D}_\alpha^z$) is equivalent to PECO. Nevertheless, it's extremely hard to figure out the exact boundary-forming data points for a complex DeDA($\mathcal{S}_\alpha^{\forall}$). Next subsection will investigate the relation between $z$ (i.e., the size of data set $\mathcal{D}^z_{\alpha}$) and the probability that the DeDA($\mathcal{D}^z_{\alpha}$) is equivalent to PECO without knowing the boundary-forming data points in advance.

\subsection{The Accuracy of Approximating PECO with DeDA($\mathcal{D}^z_{\alpha}$)} Note that DeDA($\mathcal{D}^z_{\alpha}$) is a finite deterministic program.
Following Definition \ref{def:adpfs}, we have the following definition.
\begin{definition}[on the boundary-forming data points of an optimal solution]\label{def:adpos}
The data points in the set $\mathcal{B}_{\Box}^{\rm OS}=\{\xi_{\rm d}^{(j)},\,j=1,\ldots,B_{\Box}^{\rm OS} \}$ are said to be the boundary-forming data points of $x^*_{\rm D}(\Box)$ if
\vspace{-6pt}
\begin{equation}
\vspace{-6pt}
  \mathcal{B}_{\Box}^{\rm OS}=\left\{\forall \xi_{\rm d}^{(k)} \in \mathcal{B}_{\Box}^{\rm FS} |  g_i(x^*_{\rm D}(\Box),\xi_{\rm d}^{(k)}) = 0,\;\exists i \in \{1,\ldots,m\} \right\}.
\end{equation}
\end{definition}
This definition implies that $\mathcal{B}_{\Box}^{\rm OS} \subseteq \mathcal{B}_{\Box}^{\rm FS}$ and contains all boundary-forming data points at $\mathcal{X}^*_{\rm D}(\Box)$ that each of them contributes at least one binding/active constraint to $x^*_{\rm D}(\Box)$. 
Recall that $x^*_{\rm P}$ denotes the optimal solution of the original logic model PECO (\ref{PCCSP}) and $\mathcal{S}_\alpha^{\forall}$ is the particle representation of the uncertainty set $\Upxi_\alpha$, and let $\mathcal{B}_{\forall,\alpha}^{\rm OS}$ denote the boundary-forming data set of $x^*_{\rm D}(\mathcal{S}_\alpha^{\forall})$, we have the following proposition.

\begin{proposition} \label{pro:adpos2}
    $x^*_{\rm D}(\mathcal{D}^z_{\alpha})=x^*_{\rm P}$ if $\mathcal{B}_{\forall,\alpha}^{\rm OS} 
\subseteq \mathcal{D}_\alpha^z$.
\end{proposition}
The proofs of proposition \ref{pro:adpos2} can be found in Appendix B.3. It is worth noting that $\mathcal{B}_{\forall,\alpha}^{\rm OS} 
\not \subseteq \mathcal{D}_\alpha^z$ does not necessarily mean $x^*_{\rm D}(\mathcal{D}^z_{\alpha}) \neq x^*_{\rm P}$ since there may exist other data points which are active to $x^*_{\rm D}(\mathcal{S}_\alpha^{\forall})$ and can restrict the optimal solution from changing when the corresponding data points in $\mathcal{B}_{\forall,\alpha}^{\rm OS}$ are removed. Taking Figure \ref{fig:activedata} for example, constraints 3 and 4 are the boundary-forming constraints for the optimal solution. When constraint 3 is removed, the optimal solution will not change, since active constraint 5 can provide restriction. In other words, Proposition \ref{pro:adpos2} is just a sufficient condition. Let $\varrho$ denote the probability that $x^*_{\rm D}(\mathcal{D}^z_{\alpha})$ is optimal to PECO, i.e., $\varrho=\mathbb{P}[x^*_{\rm D}(\mathcal{D}^z_{\alpha})=x^*_{\rm P}]$, and recall that $D$, $D_\alpha$, $z$, and $B_{\forall,\alpha}^{\rm OS}$ denote the numbers of data points in sets $\mathcal{D}$, $\mathcal{D}_\alpha$, $\mathcal{D}_\alpha^z$, and $\mathcal{B}_{\forall,\alpha}^{\rm OS}$ respectively, we have the following theorem on the lower bound of $\varrho$ for a given $z$ under Assumption 1. 

\textbf{Assumption 1}. $\mathcal{D}_\alpha$ \textit{contains all the boundary-forming data points of} $x^*_{\rm D}(\mathcal{S}_\alpha^{\forall})$, that is $\mathcal{B}_{\forall,\alpha}^{\rm OS} 
\subseteq \mathcal{D}_\alpha$.
\begin{theorem}\label{thm:varrho}
    Let $\bar B_{\forall,\alpha}^{\rm OS}$ be an upper bound of $B_{\forall,\alpha}^{\rm OS}$, a lower bound of $\varrho$ is 
    \begin{equation} \label{varrho}
       \underline{\varrho}(z)= 1+\sum_{k=1}^n (-1)^k\frac{\binom{\bar B_{\forall,\alpha}^{\rm OS}}{k} \binom{D_\alpha-k \alpha D}{z}}{\binom{D_\alpha}{z}}
    \end{equation}
    under Assumption 1.
\end{theorem}
 While the proof can be found in Appendix B.4, the estimation of $\bar B_{\forall,\alpha}^{\rm OS}$ and the rationality of Assumption 1 are discussed in subsection 3.4. 

 \subsection{Obtaining $\mathcal{D}_\alpha$: Defining the joint probability of a scenario based on maximum likelihood and optimal bandwidth using historical data}
According to Definition \ref{def:srd}, data points are considered i.i.d. outcomes of $\xi$, which implies that each data point in a data set shares the same probability. However, they may belong to different scenarios (see an illustrative example in Subsection 2.1). The objective of this subsection is to determine the probabilities of scenarios with only a set of historical data of $\xi$ given. For this purpose, we have the following definitions on the joint probability of a scenario through historical data points.
\begin{definition}[for integer $\xi$] \label{def:jointprob1}
Let $D_j$ be the number of data points in data set $\mathcal{D}$ that equal to scenario $\xi^{(j)}_{\rm s}$, the joint probability of $\xi^{(j)}_{\rm s}$, i.e., $P_j=\mathbb{P}[\xi=\xi^{(j)}_{\rm s}]$, is defined as the solution of the following optimization problem:
\begin{equation} \label{jointprob1}
    \max_{0\le P_j \le 1}\;\mathcal{L}(P_j)=\binom{D}{D_j}P_j^{D_j}(1-P_j)^{(D-D_j)},
\end{equation}
where $\mathcal{L}$ is the likelihood function and $\binom{\cdot}{\cdot}$ is the binomial coefficient.
\end{definition}
In the example in Section 2.1, the joint probability of $\xi^{(j)}_{\rm s}$ is defined as $D_j/D$, which is a simpler and more straightforward definition. When $\xi$ is continuous, the joint probability of a scenario is redefined as follows.
\begin{definition}[for continuous $\xi$] \label{def:jointprob2}
Let $\mathcal{D}(\xi^{(j)}_{\rm s},\zeta)$ be the set of data points in the $\zeta$-vicinity of scenario $\xi^{(j)}_{\rm s}$ in data set $\mathcal{D}$, i.e., $\mathcal{D}(\xi^{(j)}_{\rm s},\zeta)=\{\forall \xi^{(k)}_{\rm d} \in \mathcal{D} \mid \|\xi^{(k)}_{\rm d}-\xi^{(j)}_{\rm s}\| \le \zeta \}$ where $\zeta$ is a small positive scalar, and $D_j^\zeta=|\mathcal{D}(\xi^{(j)}_{\rm s},\zeta)|$, the joint probability of $\xi^{(j)}_{\rm s}$, i.e., $P_j=\mathbb{P}[\xi=\xi^{(j)}_{\rm s}]$, is defined as the solution of:
\begin{subequations} \label{jointprob2}
\begin{align}
 \min_{\zeta} \quad & \frac{1}{J}\sum_{j=1}^J\left(P_j-\mathcal{K}(\xi^{(j)}_{\rm s})\right)^2 \label{obj_djp2}   \\
 \mathrm{s.t.} \quad  &P_j=\arg\left\{ \max_{0\le P_j \le 1}\;\binom{D}{D_j^\zeta}P_j^{D_j^\zeta}(1-P_j)^{(D-D_j^\zeta)}\right\}   (j=1,\ldots, J)
\end{align} 
\end{subequations}
where $\mathcal{K}:\Upxi\rightarrow[0,1]$ is a kernel function 
and $J$ is the number of scenarios observed in $\mathcal{D}$.
\end{definition}

\begin{wrapfigure}{r}{0.446\textwidth}
\vspace{3pt}
\includegraphics[width=5.6cm]{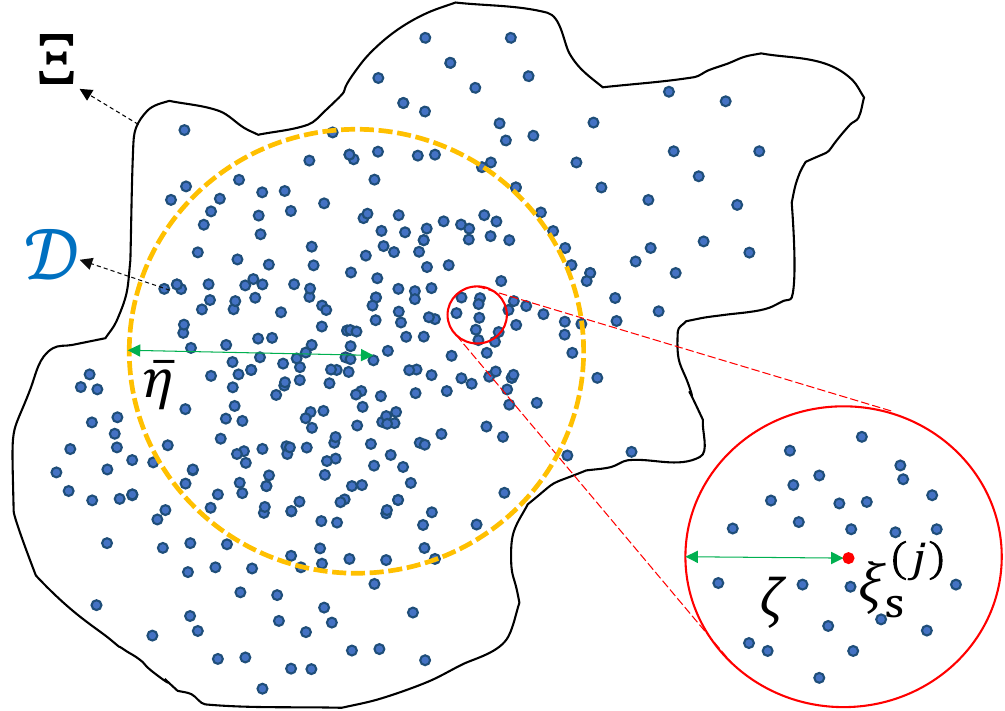}
\caption{A pictorial interpretation of $\zeta$.}
\label{fig:zeta}
\vspace{-18pt}
\end{wrapfigure}

The objective function (\ref{obj_djp2}) minimizes the deviation of the obtained probability distribution from a chosen kernel function by choosing a suitable radius $\zeta$. There are various options for the kernel function that include but are not limited to uniform, triangular, biweight, triweight, Epanechnikov (parabolic), normal, and others. $D_j^\zeta$ is a function of $\zeta$ and there may not exist a closed-form expression for this function. As shown in Figure \ref{fig:zeta}, $\zeta$ is analogous to the bandwidth of the kernel in the cases of single uncertain variable. When the uncertain parameters are mixed-integer, i.e., $\xi \in \Upxi \subset (\mathbb{Z}^{u_1},\mathbb{R}^{u_2})$ ($u_1+u_2=u$), let $\mathcal{D} |_\mathbb{Z}=\{\xi^{(k)}_{\rm d,\mathbb{Z}},\;k=1,\ldots,D\}$ be the projection of $\mathcal{D}$ in the $\mathbb{Z}^{u_1}$-space and $\mathcal{I}$ be the index set of the elements in $\mathfrak{U}[\mathcal{D}|_\mathbb{Z}]$, we divide $\mathcal{D}$ into $I=|\mathcal{I}|$ subsets, i.e., $\mathcal{D}_i$ ($i \in \mathcal{I}$), making the elements in each of these subsets have the same integer part. If scenario $\xi^{(j)}_{\rm s}$ has the same integer part as the data points in $\mathcal{D}_i$, we denote it as $\xi^{(j)}_{\rm s}=(\xi^{(i)}_{\rm s,\mathbb{Z}},\xi^{(i,k)}_{\rm s,\mathbb{R}})$ and, then have the following definition.
\vspace{-12pt}
\begin{definition}[for mixed-integer $\xi$] \label{def:jointprob3}
Let $\mathcal{D}_i(\xi^{(j)}_{\rm s},\zeta)$ be the set of data points in the $\zeta$-vicinity of scenario $\xi^{(i,k)}_{\rm s,\mathbb{R}}$ in data set $\mathcal{D}_i$, i.e., $\mathcal{D}_i(\xi^{(j)}_{\rm s},\zeta)=\{\forall \xi^{(i,l)}_{\rm d} \in \mathcal{D}_i \mid \|\xi^{(i,l)}_{\rm d,\mathbb{R}}-\xi^{(i,k)}_{\rm s,\mathbb{R}}\| \le \zeta \}$, and $D_j^\zeta=|\mathcal{D}_i(\xi^{(j)}_{\rm s},\zeta)|$, the joint probability of $\xi^{(j)}_{\rm s}$, i.e., $P_j=\mathbb{P}[\xi=\xi^{(j)}_{\rm s}]$, is defined as the solution of (\ref{jointprob2}).
\end{definition}
Data set $D_{\alpha}$ can be obtained by picking all data points in $\mathcal{D}$ that belong to scenarios whose joint probabilities are not less than $\alpha$ as defined in definitions \ref{def:jointprob1}-\ref{def:jointprob3}.

\subsection{Discussions and Summary} In Section 3, we considered DeDA($\mathcal{D}_\alpha^z$) as the deterministic approximation PECO (\ref{PCCSP}) and revealed some of its properties. This subsection provides a further discussion and a summary.
\subsubsection{A proper estimation of $\bar B_{\forall,\alpha}^{\rm OS}$}
It's generally very hard to obtain the exact number of boundary-forming data points for the optimal solution of a D-DA. Therefore, it is necessary to obtain a proper estimate of $\bar B_{\forall,\alpha}^{\rm OS}$ before the assertion of Theorem \ref{thm:varrho} can be used to determine $z$ needed for a desired $\rho$. We start discussing how to obtain a proper upper bound of the number of boundary-forming data points for a specific case from the following proposition.
\begin{proposition}\label{pro:adpos1}
    For a special case of DeDA (\ref{DDA}), where $f(x)$ is linear and $g(x,\xi)$ is convex on $x$, the number of boundary-forming data points of its optimal solution is not more than $n$, where $n$ is the size of the vector of decision variables $x$.
\end{proposition}
While the proof is provided in Appendix B.5, this proposition asserts that, for DeDAs which are convexly-constrained linear programs, $\bar B_\Box^{\rm OS}=n$ is a valid upper bound of $B_\Box^{\rm OS}$. Although this assertion may not be directly applicable to a more general case, it can provide a proper estimation. Because, from the proof of Proposition \ref{pro:adpos1}, we know that $n$ is in fact a very loose upper bound of $B_\Box^{\rm OS}$ for convexly constrained linear programs. To be specific, ``$n \ge$ the number of support constraints $\ge$ the number of boundary-forming constraints $\ge$ the number of boundary-forming data points." Therefore, $n$ is a valid upper bound of $B_\Box^{\rm OS}$ (including $B_{\forall,\alpha}^{\rm OS}$) for a much wider range of cases than the convexly constrained linear programs although it cannot be guaranteed for all cases.
\subsubsection{Rationality of the data-embedded deterministication}
We first discuss the rationality of the direct use of historical data as input for determinisfying a PECO problem. In engineering practice, the construction of the probability triple, i.e., ($\Upxi$, $\Upxi$, $P$), is generally based on engineers' experience, which comes mainly from the analysis on historical data and/or observations together with some assumptions. While inaccuracy can occur in the process of constructing the probability triple, such an issue does not exist in the proposed data-embedded method since it does not need the process of constructing probability triples. 
Second, we'd like to discuss the rationality of Assumption 1. Normally, engineers care about the outcomes of $\xi$ which occurred in history, especially those that frequently appeared. If $\mathcal{D}$ is properly measured in history, $\mathcal{D}_\alpha$ should contain all high-probability outcomes of $\xi$ that engineers are interested in and, of course, contain all boundary-forming data points. In other words, we can consider that, if a scenario of $\xi$ never occurred in history, it's neither of interest nor a boundary-forming data point. In light of this, we have a reasonable hypothesis that the boundary-forming data points defined in Definitions \ref{def:adpfs} and \ref{def:adpos} are generally some frequently appeared data points in history.

\section{The Proposed Solution Paradigm for PECO: Strategic Data Selection}

If DeDA($\mathcal{D}_\alpha^z$), that was investigated in Section 3, is computationally intractable, the proposed solution paradigm will use DeDA($\mathcal{D}_\alpha^\eta$) as the deterministic approximation of PECO instead, where $\mathcal{D}_\alpha^{\eta} \subseteq \mathcal{D}_\alpha^z$. In other words, a subset $\mathcal{D}_\alpha^{\eta}$ of $\mathcal{D}_\alpha^z$ is embedded into DeDA (\ref{DDA}). Since $\mathcal{D}_\alpha^{\eta} \subseteq \mathcal{D}_\alpha^z \subseteq \mathcal{D}_\alpha$, it's straightforward to know that, the feasible sets satisfy:
\begin{equation} \label{FSrelations}
    \mathcal{X}_{\rm D}(\mathcal{D}_\alpha^{\eta}) \supseteq \mathcal{X}_{\rm D}(\mathcal{D}_\alpha^z) \supseteq \mathcal{X}_{\rm D}(\mathcal{D}_\alpha)
\end{equation}
and, hence, $f(x^*_{\rm D}(\mathcal{D}_\alpha^{\eta})) \le f(x^*_{\rm D}(\mathcal{D}_\alpha^z)) \le f(x^*_{\rm D}(\mathcal{D}_\alpha))$ if both $f$ and $g$ are convex in $x$. Proposition \ref{pro:adpos2} indicates that $x^*_{\rm D}(\mathcal{D}_\alpha)=x^*_{\rm P}$ under Assumption 1, and hence $f(x^*_{\rm D}(\mathcal{D}_\alpha^{\eta})) \le f(x^*_{\rm D}(\mathcal{D}_\alpha^z)) \le f(x^*_{\rm P})$. While the probabilistic gap between $x^*_{\rm D}(\mathcal{D}_\alpha^z)$ and $x^*_{\rm P}$, i.e., $\mathbb{P}[x^*_{\rm D}(\mathcal{D}_\alpha^z) \neq x^*_{\rm P}]$, was discussed in the previous section, this section aims to develop a methodology of strategic data selection (SDS) for obtaining the data set $\mathcal{D}_\alpha^{\eta}$ which possesses two features: 1) the number of elements $z_\eta =|\mathcal{D}_\alpha^\eta| \ll z$, and 2) the Euclidean gap between $x^*_{\rm D}(\mathcal{D}_\alpha^{\eta})$ and $x^*_{\rm D}(\mathcal{D}_\alpha^z)$, i.e., $\|x^*_{\rm D}(\mathcal{D}_\alpha^{\eta})-x^*_{\rm D}(\mathcal{D}^z_{\alpha})\|$, is sufficiently small such that $x^*_{\rm D}(\mathcal{D}_\alpha^{\eta})$ can also be considered a good estimation to $x^*_{\rm P}$. The SDS algorithms for continuous and integer/mixed-integer cases of $\xi$ are presented in Subsections 4.1 and 4.3, respectively.

\subsection{The SDS algorithm for continuous $\xi$}
When the uncertain parameters are continuous, i.e., $\xi \in \Upxi \subset \mathbb{R}^u$, we developed an SDS algorithm as detailed in Algorithm \ref{alg:sds}, where $\eta$ is a key parameter. Then, we have the following theorem on the relation between $z_\eta$ and $\eta$.

\begin{theorem}[on $z_\eta$ obtained by Algorithm \ref{alg:sds}] \label{pro:continuous}
  When $\xi \in \Upxi \subset \mathbb{R}^u$, we have
   \begin{equation} 
    z_\eta \begin{cases} =z &  \eta=0  \\\le \Bar{z}_\eta & 0< \eta < \bar \eta \\ =1 & \eta \ge \bar \eta \end{cases}. \nonumber
    \end{equation}
When $0< \eta < \bar \eta$, $z_\eta$ is a parameterized random variable whose expectation and upper bound satisfy the following relations respectively
    \begin{equation}  \label{pro_con_2}
    \frac{\partial (\mathbb{E}[z_\eta])}{\partial \eta}   \le 0,\; \text{and}\;\; \Bar{z}_\eta \propto \frac{1}{\eta^u}.
    \end{equation}
\end{theorem}
The proof of theorem \ref{pro:continuous} is provided in Appendix C.1. The second expression in (\ref{pro_con_2}) indicates that, when $\eta$ is small, $z_\eta$ drops exponentially as $\eta$ increases.

Given that a data point $\xi_{\rm d}^{(k)}$ contributes $m$ constraints to DeDA (\ref{DDA}) since $g:\mathbb{R}^{n+u}  \rightarrow \mathbb{R}^m$ in the PEC (\ref{PCC}), there are $z \times m$ inequality constraints in DeDA($\mathcal{D}_\alpha^z$). Hence, let $\xi_{\rm d,B}^{(k)}$ be a boundary-forming data point of the optimal solution $x^*_{\rm D}(\mathcal{D}_\alpha^z)$ (as defined in Definition \ref{def:adpos}), it contributes $m$ constraints and, at least, one boundary-forming constraint to DeDA($\mathcal{D}_\alpha^z$). Among these $m$ constraints, if $\Tilde{g}_i(x,\xi_{\rm d,B}^{(k)}) \le 0$ ($i \in 1,\ldots,m$) is a boundary-forming constraint, the equal sign holds when $x$ is substituted by $x^*_{\rm D}(\mathcal{D}^z_{\alpha})$, that is $\Tilde{g}_i(x^*_{\rm D}(\mathcal{D}^z_{\alpha}),\xi_{\rm d,B}^{(k)}) = 0$. Suppose there are $B_c$ boundary-forming constraints for $x^*_{\rm D}(\mathcal{D}^z_{\alpha})$, we have:
\begin{equation} \label{active}
    \Tilde{g}(x^*_{\rm D}(\mathcal{D}^z_{\alpha}),\xi_{\rm d,B})= \begin{cases} \Tilde{g}_1(x^*_{\rm D}(\mathcal{D}^z_{\alpha}),\xi_{\rm d,B}^{(1)})=0\\ \vdots \\ \Tilde{g}_{B_c}(x^*_{\rm D}(\mathcal{D}^z_{\alpha}),\xi_{\rm d,B}^{(B(z))})=0 \end{cases}\text{, where}\;\Tilde{g}(x,\Tilde{\xi})= \begin{cases} \Tilde{g}_1(x,\xi_1)\\ \vdots \\ \Tilde{g}_{B_c}(x,\xi_{B_c}), \end{cases} 
    \end{equation}
$\Tilde{g}:\mathbb{R}^{n+(u\times B_c)} \rightarrow \mathbb{R}^{B_c}$ ($\Tilde{g}_i:\mathbb{R}^{n+u} \rightarrow \mathbb{R}$ and $\xi_i \in \mathbb{R}^u$ ($i=1,\ldots,B_c$)), $\Tilde{\xi}=[\xi_1^{\rm T},\ldots, \xi_{B_c}^{\rm T}]^{\rm T}$ $\in \mathbb{R}^{(u\times B_c) \times 1}$ and $\xi_{\rm d,B}=[\xi_{\rm d,B}^{(1)\rm T},\ldots, \xi_{\rm d,B}^{(B(z))\rm T}]^{\rm T} \in \mathbb{R}^{(u\times B_c) \times 1}$ which implies that the $B_c$ boundary-forming constraints are contributed by $B(z)$ boundary-forming data points (note that $B_c \ge B(z)$). In other words, it's possible that two or more elements in $\xi_{\rm d,B}$ are identical, which means that this data point (say $\xi_{\rm d,B}^{(i)}$) contributes more than one boundary-forming constraints. Constraints (\ref{active}) are the $B_c$ constraints out of $z \times m$ whose equal sign holds for the optimal solution $x^*_{\rm D}(\mathcal{D}^z_{\alpha})$. According to the discussion in Subsubsection 3.4.1, it's reasonable to assume that $B_c \le n$. Then, we have the following proposition.
\begin{algorithm}[t]
\caption{of strategic data selection for continous $\xi$}
\label{alg:sds}
\begin{algorithmic}
\STATE{Given a data set $\mathcal{D}_\alpha^z$:}
\begin{enumerate}[leftmargin=*]
  \item randomly select a data point $\xi^{(0)}_{\rm d}$ and set $i=1$;
  \item randomly select another data point $\xi^{(i)}_{\rm d}$ in $\mathcal{D}^z_\alpha$ which satisfies
  \begin{equation} \label{agr2}
      \|\xi^{(i)}_{\rm d} -\xi^{(j)}_{\rm d} \| \ge 2\eta,\; \forall j=0,\ldots, i-1
  \end{equation}
  where $0 < \eta \le \bar \eta$;
  \item save $\xi^{(i)}_{\rm d}$ to a new set $\mathcal{D}_\alpha^{\eta}$ and discard all data points in the $\eta$-vicinity of $\xi^{(i)}_{\rm d}$ (including $\xi^{(i)}_{\rm d}$) from $\mathcal{D}^z_\alpha$;
  \item stop and report set $\mathcal{D}_\alpha^{\eta}$ if there is no data point $\xi^{(i)}_{\rm d}$ in $\mathcal{D}^z_\alpha$ that satisfies (\ref{agr2}), otherwise go to step 5;
  \item set $i=i+1$ and repeat steps 2-4.
\end{enumerate} 
\STATE{Note: $\bar \eta$ is the minimum $\eta$ which results in that there do not exist two data points in $\mathcal{D}^z_\alpha$ which satisfy condition (\ref{agr2}) (see an illustrative example in Figure \ref{fig:zeta}).} 
\end{algorithmic}
\end{algorithm}
\begin{proposition}[Implicit function theorem]\label{pro:ift}
     When $B_c = n$, there exists a vector-valued function $x=h(\Tilde{\xi})$ ($h:\mathbb{R}^{(u \times n)} \rightarrow \mathbb{R}^n$) which is equivalent to $\Tilde{g}(x,\Tilde{\xi})=0$ ($\Tilde{g}:\mathbb{R}^{n+(u \times n)} \rightarrow \mathbb{R}^n$) in the vicinity of ($x^*_{\rm D}(\mathcal{D}^z_{\alpha}),\xi_{\rm d,B}$) if $\Tilde{g}$ is continuously differentiable and the Jacobian matrix of $\Tilde{g}$ with respect to $x$ evaluated at $x^*_{\rm D}(\mathcal{D}^z_{\alpha})$ is invertible.
\end{proposition}
Proposition \ref{pro:ift} is actually the Implicit Function Theorem whose proof can be found in \cite{krantz2002implicit}.  Algorithm \ref{alg:sds} guarantees that, for any boundary-forming data point $\xi_{\rm d,B}^{(i)}$ ($i=1,\ldots,B(z)$) of $x^*_{\rm D}(\mathcal{D}^z_{\alpha})$, $\mathcal{D}_\alpha^{\eta}$ contains a data point $\acute{\xi}_{\rm d}^{(i)}$ which satisfies $\|\acute{\xi}_{\rm d}^{(i)}- \xi_{\rm d,B}^{(i)}\|\le 2\eta$. Recall that $x^*_{\rm D}(\mathcal{D}_\alpha^{\eta})$ denotes the optimal solution of DeDA($\mathcal{D}_\alpha^{\eta}$) and let $\varphi=1-\frac{\|x^*_{\rm D}(\mathcal{D}_\alpha^{\eta})-x^*_{\rm D}(\mathcal{D}^z_{\alpha})\|}{\| x^*_{\rm D}(\mathcal{D}^z_{\alpha})\|}$ denote the accuracy\footnote{If needed, $\varphi=1-\frac{\|f(x^*_{\rm D}(\mathcal{D}_\alpha^{\eta}))-f(x^*_{\rm D}(\mathcal{D}^z_{\alpha}))\|}{\| f(x^*_{\rm D}(\mathcal{D}^z_{\alpha}))\|}$ can be used to represent this accuracy.} of approximating DeDA($\mathcal{D}_\alpha^z$) with DeDA($\mathcal{D}_\alpha^{\eta}$), we have $\varphi=1$ when DeDA($\mathcal{D}_\alpha^{\eta}$) is equivalent to DeDA($\mathcal{D}_\alpha^z$). It's straightforward to know that $\varphi$ is partially determined by $\eta$ for a given $\mathcal{D}_\alpha^z$. Denoting $\acute{\xi}_{\rm d}^{(i)}$ as the closest data point to $\xi_{\rm d,B}^{(i)}$ ($\forall i=1,\dots,B_c$) in $\mathcal{D}_\alpha^{\eta}$, we have $\acute{\xi}_{\rm d}=[\acute{\xi}_{\rm d}^{(1)\rm T},\ldots, \acute{\xi}_{\rm d}^{(B_c)\rm T}]^{\rm T}$ $\in \mathbb{R}^{(u\times B_c) \times 1}$. Further let  $\xi_{\rm d,B^\prime}^{(i)}$ ($i=1,\ldots,B(\eta)$) denote the boundary-forming data points of $x^*_{\rm D}(\mathcal{D}_\alpha^{\eta})$, we have the following theorem.
\begin{theorem}[on a lower bound of $\varphi$ under Algorithm \ref{alg:sds}] \label{thm:rho2}
     When $\eta$ is small, a lower bound of $\varphi$ is given as (\ref{thm4.3}) under the conditions in Proposition \ref{pro:ift}, 
     \begin{equation} \label{thm4.3}
         \underline{\varphi}(\eta)=1-\frac{2\sqrt{n}\|H(\xi_{\rm d,B})\|\eta+\frac{2^2n}{2!}\|H^\prime(\xi_{\rm d,B})\|\eta^2+\frac{2^3\sqrt{n^3}}{3!}\|H^{\prime\prime}(\xi_{\rm d,B})\|\eta^3+\cdots}{\| x^*_{\rm D}(\mathcal{D}_\alpha^z)\|} 
     \end{equation}
     where $H(\xi_{\rm d,B})$ and $H^\prime(\xi_{\rm d,B})$ are the Jacobian matrix and Hessian tensor, respectively, of $h(\xi)$ evaluated at $\xi_{\rm d,B}$, and $H^{\prime\prime}(\xi_{\rm d,B})$ is a tensor which is a higher-dimensional generalization of a matrix and contains the third mixed partial derivatives, 
     if the problem satisfies any of the following conditions:
\begin{enumerate}[leftmargin=*]
    \item $\acute{\xi}_{\rm d}^{(i)}$ ($i=1,\dots,B_c$) are exactly the boundary-forming data points of $x^*_{\rm D}(\mathcal{D}_\alpha^{\eta})$, i.e., $\{\xi_{\rm d,B^\prime}^{(i)}\,(i=1,\ldots,B(\eta))\}=\{\acute{\xi}_{\rm d}^{(i)}\,(i=1,\dots,B_c)\}$.
    \item The function $g$ is convex in $x$ and the objective function $f$ is convex and radially non-increasing in the vicinity of $x^*_{\rm D}(\mathcal{D}_\alpha^z)$ alone any direction $d\in \mathbb{R}^n$ that satisfies
\begin{equation} \label{outward}
    \nabla\Tilde{g}_i(x^*_{\rm D}(\mathcal{D}_\alpha^z),\xi_{\rm d,B}^{(i)})^{\rm T}d 
\ge 0,\, \forall i=1,\dots,B_c  
\end{equation}
where $\nabla$ denotes the derivative of a function with respect to $x$.   
\end{enumerate}
\end{theorem}
The proof of Theorem \ref{thm:rho2} is provided in Appendix C.2. Given the continuous differentiability of $g$ and $f$, condition 1 is easy to satisfy when $\eta$ is small. In condition 2, a direction that satisfies (\ref{outward}) is an ``outward" pointing direction for the boundary-forming constraints at $x^*_{\rm D}(\mathcal{D}_\alpha^z)$. A discussion on the benefits of SDS is provided in next subsection. 
\subsection{Discussion on the benefits of SDS} First, it's worth noting that the computational complexity of the DeDAs (\ref{DDA}) is related to their numbers of embedded data points, i.e., the size of $\Box$. For simple cases, such as linear cases, the computational complexity is approximately proportional to the number of embedded data points. However, for non-convex cases, the computational burden generally grows exponentially as the number of embedded data points increases. Second, from theorems \ref{pro:continuous} and \ref{thm:rho2}, we know that both $z_\eta$, which is related to computational complexity, and $\varphi$, which is related to precision, decrease as $\eta$ increases. To be specific, Theorem \ref{thm:rho2} implies (\ref{barh2}), which means that, when $\eta$ is small, $\underline{\varphi}(\eta)$ is close to being linear to $\eta$. By comparing the expressions of $\bar z_\eta$ and $\underline \varphi$, we realize that $\bar z_\eta$ drops significantly while $\underline \varphi$ decreases at a much slower rate. Therefore, the main benefit of SDS is that it significantly reduces the computational burden of the deterministic approximation of PECO, i.e., DeDA($\mathcal{D}_\alpha^{\eta}$), with a tradeoff of slightly decrease in accuracy.

\begin{equation}  
    \underline{\varphi}(\eta) \approx \begin{cases} 1-\frac{2\sqrt{n}\|H(\xi_{\rm d,B})\|\eta}{\| x^*_{\rm D}(\mathcal{D}_\alpha^z)\|} & \underline \eta \le \eta \le \eta_1 \\ 1-\frac{2\sqrt{n}\|H(\xi_{\rm d,B})\|\eta+\frac{2^2n}{2!}\|H^\prime(\xi_{\rm d,B})\|\eta^2}{\| x^*_{\rm D}(\mathcal{D}_\alpha^z)\|} & \eta_1 < \eta \le \eta_2  \\ 1-\frac{2\sqrt{n}\|H(\xi_{\rm d,B})\|\eta+\frac{2^2n}{2!}\|H^\prime(\xi_{\rm d,B})\|\eta^2+\frac{2^3\sqrt{n^3}}{3!}\|H^{\prime\prime}(\xi_{\rm d,B})\|\eta^3}{\| x^*_{\rm D}(\mathcal{D}_\alpha^z)\|} & \eta_2 < \eta \le \eta_3 \\  \vdots & \vdots \\1-\frac{2\sqrt{n}\|H(\xi_{\rm d,B})\|\eta+\frac{2^2n}{2!}\|H^\prime(\xi_{\rm d,B})\|\eta^2+\frac{2^3\sqrt{n^3}}{3!}\|H^{\prime\prime}(\xi_{\rm d,B})\|\eta^3+\cdots}{\| x^*_{\rm D}(\mathcal{D}_\alpha^z)\|}  & \eta_j < \eta \le \bar \eta \end{cases}, \label{barh2}
    \end{equation}

\subsection{The SDS algorithm for integer/mixed-integer $\xi$}
When the uncertain parameters are mixed-integer, i.e., $\xi \in \Upxi \subset (\mathbb{Z}^{u_1},\mathbb{R}^{u_2})$ ($u_1+u_2=u$), let $\mathcal{D}_\alpha^z |_\mathbb{Z}=\{\xi^{(k)}_{\rm d,\mathbb{Z}},\;k=1,\ldots,z\}$ be the projection of $\mathcal{D}_\alpha^z$ in the $\mathbb{Z}^{u_1}$-space, and $\mathcal{I}$ be the index set of the elements in $\mathfrak{U}[\mathcal{D}_\alpha^z|_\mathbb{Z}]$, we divide $\mathcal{D}_\alpha^z$ into $I=|\mathcal{I}|$ subsets, i.e., $\mathcal{D}^z_{\alpha,i}$ ($i \in \mathcal{I}$), making the elements in each of these subsets have the same integer part. Namely, if we denote all data points in $\mathcal{D}^z_{\alpha,i}$ as $\xi^{(i,j)}_{\rm d}=(\xi^{(i)}_{\rm d,\mathbb{Z}},\xi^{(i,j)}_{\rm d,\mathbb{R}})$ ($j=1,\ldots,z_i$), they share the same integer part, i.e., $\xi^{(i)}_{\rm d,\mathbb{Z}}$, where $z_i$ is the number of data points in $\mathcal{D}^z_{\alpha,i}$. Based on these notations, Algorithm \ref{alg:mix} of SDS is developed for the cases where $\xi$ is integer/mixed-integer. Then, we have the following corollary of Theorem \ref{pro:continuous}.
\begin{algorithm}
\caption{of strategic data selection for cases of integer/mixed-integer $\xi$}
\label{alg:mix}
\begin{algorithmic}
\STATE{Given a data set $\mathcal{D}^z_\alpha=\bigcup_{i \in \mathcal{I}}\mathcal{D}^z_{\alpha,i}$:}
\begin{enumerate}[leftmargin=*]
  \item set $i=1$
 \item randomly select a data point $\xi^{(i,0)}_{\rm d}$ in $\mathcal{D}^z_{\alpha,i}$ and set $j=1$;
     \item randomly select another data point $\xi^{(i,j)}_{\rm d}$ in $\mathcal{D}^z_{\alpha,i}$ which satisfies
  \begin{equation} \label{agr3}
      \|\xi_{\rm d,\mathbb{R}}^{(i,j)} -\xi_{\rm d,\mathbb{R}}^{(i,k)} \| \ge 2\eta_i,\; \forall k=0,\ldots, j-1;
  \end{equation}
  \item save $\xi^{(i,j)}_{\rm d}$ to a new set $\mathcal{D}^\eta_{\alpha}$ and discard all data points in the $\eta_i$-vicinity of $\xi^{(i,j)}_{\rm d}$ (including $\xi^{(i,j)}_{\rm d}$) from $\mathcal{D}^z_{\alpha,i}$;
  \item go to step 7 if there is no data points in $\mathcal{D}^z_{\alpha,i}$ that satisfies (\ref{agr3}), otherwise go to step 6;
  \item set $j=j+1$ and repeat steps 3-5.
  \item stop and report set $\mathcal{D}^\eta_{\alpha}$ if $i=|\mathcal{I}|$, otherwise go to step 8;
  \item set $i=i+1$ and repeat steps 2-7.
\end{enumerate} 
\STATE{Note: $\bar \eta_i$ is the minimum $\eta_i$ which results in that there do not exist two data points in $\mathcal{D}^z_{\alpha,i}$ which satisfy condition (\ref{agr3}).}
\end{algorithmic}
\end{algorithm}

\begin{corollary}[on $z_\eta$ obtained by Algorithm \ref{alg:mix}] \label{cor:Neta}
  When $\xi \in \Upxi \subset (\mathbb{Z}^{u_1},\mathbb{R}^{u_2})$ ($u_1+u_2=u$), we have
   \begin{equation}  \label{cor_con_1}
    z_\eta \begin{cases} =z &  \eta_i =0 \\\le \Bar{z}_\eta & 0< \eta_i < \bar \eta_i \\ =I & \eta_i \ge \bar \eta_i \end{cases}, \; i \in \mathcal{I}.
    \end{equation}
where $I=|\mathcal{I}|=|\mathfrak{U}[\mathcal{D}_\alpha^z |_\mathbb{Z}]|$. When $0< \eta_i < \bar \eta_i$ ($\forall i \in \mathcal{I}$), $z_\eta$ is a parameterized random variable whose expectation and upper bound satisfy the following relations respectively
    \begin{equation}  \label{cor_con_2}
    \frac{\partial (\mathbb{E}[z_\eta])}{\partial \eta_i}   \le 0,\; \; \Bar{z}_\eta=\sum_i \Bar{z}_{\eta,i}\;\text{and}\;\Bar{z}_{\eta,i} \propto \frac{1}{\eta_i^{r_2}}\;\;(i \in \mathcal{I}).
    \end{equation}
\end{corollary}
The proof of Corollary \ref{cor:Neta} is in Appendix C.3. 

For the cases where $\xi$ is mixed-integer, we also denote the boundary-forming constraints of $x^*_{\rm D}(\mathcal{D}_\alpha^z)$ in DeDA($\mathcal{D}_\alpha^z$) as $\Tilde{g}(x,\xi_{\rm d,B}) \le 0$ ($\Tilde{g}:\mathbb{R}^n\times \Upxi \rightarrow \mathbb{R}^{B_c}$), which are contributed by $B(z)$ boundary-forming data points. Let $\xi_{\rm d,B}^{(k)}$ ($k \in 1,\ldots,B(z)$) denote a boundary-forming data point of $x^*_{\rm D}(\mathcal{D}_\alpha^z)$ and assume that it belongs to subset $\mathcal{D}^z_{\alpha,i}$ ($i \in \mathcal{I}$). When $\eta_i$ ($\forall i \in \mathcal{I}$) are small, algorithm \ref{alg:mix} guarantees that $\mathcal{D}_\alpha^{\eta}$ contains a data point $\acute{\xi}_{\rm d}^{(k)}$ which satisfies $\acute{\xi}_{\rm d,\mathbb{Z}}^{(k)}= \xi_{\rm d,B,\mathbb{Z}}^{(k)}$ and $\|\acute{\xi}_{\rm d,\mathbb{R}}^{(k)}- \xi_{\rm d,B,\mathbb{R}}^{(k)}\|\le 2\eta_i$. Given that, for mixed-integer $\xi$, the data points in a small vicinity of $\xi_{\rm d}^{(k)}$ have the same integer part as $\xi_{\rm d}^{(k)}$, one can consider $\xi_{\mathbb{Z}} \in \mathbb{Z}^{u_1}$ fixed and reformulate $\Tilde{g}(x,\Tilde{\xi})$ as $\breve{g}(x,\Tilde{\xi}_{\mathbb{R}})$ ($\breve{g}:\mathbb{R}^{n+(u_2\times B_c)} \rightarrow \mathbb{R}^{B_c}$). According to Proposition \ref{pro:ift} and under the conditions therein, we know that there exists a vector-valued function $x=h(\Tilde{\xi}_{\mathbb{R}})$ ($h:\mathbb{R}^{(u_2 + n)} \rightarrow \mathbb{R}^n$) which is equivalent to $\breve{g}(x,\Tilde{\xi}_{\mathbb{R}})=0$ in the vicinity of ($x^*_{\rm D}(\mathcal{D}_\alpha^z),\xi_{\rm d,B}$) (Note that $B_c=n$ under the condition in Proposition \ref{pro:ift}.). Then, we have the following corollary of Theorem \ref{thm:rho2}.
\begin{corollary}[on a lower bound of $\varphi$ under Algorithm \ref{alg:mix}] \label{cor:rho2}
     When $\eta_i$ ($\forall i \in \mathcal{I}$) are small, a lower bound of $\varphi$ is given as follows:  
     \begin{equation} \label{cor_rho2}
         \underline{\varphi}(\eta_1,\ldots,\eta_{B_c})=1-\frac{\|H(\xi_{\rm d,B,\mathbb{R}})\|\hat{\eta}+\frac{1}{2!}\|H^\prime(\xi_{\rm d,B,\mathbb{R}})\|\hat{\eta}^2+\frac{1}{3!}\|H^{\prime\prime}(\xi_{\rm d,B,\mathbb{R}})\|\hat{\eta}^3+\cdots}{\| x^*_{\rm D}(\mathcal{D}_\alpha^z)\|} ,
     \end{equation}
     where $\hat{\eta}=2\sqrt{\eta_1^2+\ldots+\eta_{B_c}^2}$, under the conditions in Proposition \ref{pro:ift} if the problem satisfies any of the two conditions in Theorem \ref{thm:rho2}.
\end{corollary}
The proof of Corollary \ref{cor:rho2} can be found in Appendix C.4. 

When the uncertain parameters are pure integers, i.e., $\xi \in \Upxi \subset \mathbb{Z}^u$, step 3 in Algorithm \ref{alg:mix} is not necessary since $u_2=0$ and, consequently, $\xi_{\rm d,\mathbb{R}}=0$. Moreover, the data points in each subset $\mathcal{D}^z_{\alpha,i}$ ($i \in \mathcal{I}$) are identical. Thus, Algorithm \ref{alg:mix} picks one data point from each of these subsets, which results in that $\mathcal{D}_\alpha^{\eta}$ is the underlying set of $\mathcal{D}_\alpha^z$, i.e., $\mathcal{D}_\alpha^{\eta}=\mathfrak{U}[\mathcal{D}_\alpha^z]$, and $z_\eta=I$. It's worth noting that, generally, $I << z$ (see the illustrative example in Figure \ref{fig:realization_scenario}). Then, we have the following lemma on Algorithm \ref{alg:mix} for the cases of pure integer $\xi$.

\begin{lemma}[on SDS algorithm for discrete $\xi$] \label{pro:discrete}
     When $\xi \in \Upxi \subset \mathbb{Z}^u$, $x^*_{\rm D}(\mathcal{D}_\alpha^{\eta})=x^*_{\rm D}(\mathcal{D}_\alpha^z)$. Further, if $\mathcal{D}_\alpha^z=\mathcal{D}_\alpha$ (i.e., Algorithm \ref{alg:mix} uses $\mathcal{D}_\alpha$ as the input data set),  $x^*_{\rm D}(\mathcal{D}_\alpha^{\eta})=x^*_{\rm P}$ under Assumption 1. 
\end{lemma}
The proof of this lemma is provided in Appendix C.5. Let $\omega$ denote the probability that $x^*_{\rm D}(\mathcal{D}_\alpha^{\eta})$ is optimal to PECO, i.e., $\omega=\mathbb{P}[x^*_{\rm D}(\mathcal{D}_\alpha^{\eta})=x^*_{\rm P}]$, Lemma \ref{pro:discrete} implies that $\omega=1$ when Algorithm \ref{alg:mix} uses $\mathcal{D}_\alpha$ as input data set and Assumption 1 holds. The following proposition presents a trivial property of $\omega$ for the PECO problems with continuous/mixed-integer $\xi$, of which the proof is given in Appendix C.6.

\begin{proposition}[on $\omega$ for continuous/mixed-integer $\xi$]  \label{pro:omega}
If $\eta_1 < \eta_2$, we have $\mathbb{P}[\omega(\eta_1) > \omega(\eta_2)]\ge\mathbb{P}[\omega(\eta_1) < \omega(\eta_2)]$.
\end{proposition}

\section{Application and Numerical Experiments}
This section applies the proposed methods to a fundamental decision-making problem under uncertainty in electric power systems, i.e., the optimal power flow (OPF) with uncertain renewable energy, e.g., solar and wind power. The OPF determines the best operating levels of power generators/plants in order to meet demands given throughout a transmission/distribution network, usually with the objective of minimizing generation cost \cite{li2016convex}. The formulation of the deterministic OPF is given as 
\vspace{-6pt}
\begin{subequations} \label{OPF} 
\begin{align}
\text{\textbf{OPF}:}\quad\quad \min_{p^G,\theta} \quad & c(p^G)=\sum_{i\in \mathcal{G}} \left(c_{i,2}(p^G_i)^2+c_{i,1}p^G_i+c_{i,0}\right)  \\
 \mathrm{s.t.} \quad  & B \theta =A p^G+C p^R+d, \label{DCPF} \\
 & \left|B_{ij}(\theta_i-\theta_j)\right| \le P_{ij}^{\max},\; \forall\{i,j \} \in \mathcal{E}\\
 & p_i^{\min} \le p^G_i \le p_i^{\max}, \forall i\in \mathcal{G}
\end{align} 
\end{subequations}
where (\ref{DCPF}) is the DC power flow (DCPF) equation, and A, B, and C are $n \times g$, $n \times n$, and $n \times u$ matrices whose elements are given as follows respectively:
\vspace{-6pt}
\begin{align}
&  A_{ij}=\begin{cases}
         1, & \text{if the $j$th power generator is connected to the $i$th node} \\
          0, & \text{otherwise}
    \end{cases} \nonumber \\
    & B_{ij}=\begin{cases}
         -b_{ij}, & \{i,j\} \in \mathcal{E} \\
          \sum_{k;\{k,j\}\in \mathcal{E}}b_{kj}, & i=j \\
          0, & \text{otherwise}
    \end{cases}  \nonumber \\
    &  C_{ij}=\begin{cases}
         1, & \text{if the $j$th renewable generator is connected to the $i$th node} \\
          0, & \text{otherwise}
    \end{cases} \nonumber
\end{align}
   \vspace{-6pt}
A nomenclature is given in the Table \ref{tab:nct}.

\begin{table}[h]
\footnotesize
\centering
\caption{Nomenclature for OPFRG}
\vspace{-6pt}
\begin{tabular}{|lp{0.88\linewidth}|}
\hline
\multicolumn{2}{|l|}{\textit{A. Sets and Indices}} \\ \hline
\multicolumn{1}{|l|}{$\mathcal{E}$}         &   Set of transmission/distribution lines       \\ \hline
\multicolumn{1}{|l|}{$\mathcal{G}$}         &  Set of $g$ power generators/plants        \\ \hline
\multicolumn{1}{|l|}{$\mathcal{N}$}         &  Set of $n$ nodes/buses        \\ \hline
\multicolumn{1}{|l|}{$\mathcal{R}$}         &  Set of $u$ renewable generators        \\ \hline
\multicolumn{2}{|l|}{\textit{B. Parameters}}       \\ \hline
\multicolumn{1}{|l|}{$b_{ij}$}         &  Susceptance of the transmission line $\{i,j\} \in \mathcal{E}$ \\ \hline
\multicolumn{1}{|l|}{$c_i$}         &  Unit fuel cost of the $i$th power plant in $\$$/MWh, where $i \in \mathcal{G}$     \\ \hline
\multicolumn{1}{|l|}{$d$}         &  $n \times 1$ vector of electricity demands     \\ \hline
\multicolumn{1}{|l|}{$P_{i,j}^{\max}$}    &  Power limit on the transmission line $\{i,j\} \in \mathcal{E}$  \\ \hline
\multicolumn{1}{|l|}{$P_i^{\min}$}   &   Lower bound on power generation of $i$th generator/plant, where $i \in \mathcal{G}$   \\ \hline
\multicolumn{1}{|l|}{$P_i^{\max}$}    &  Upper bound on power generation of $i$th generator/plant, where $i \in \mathcal{G}$   \\ \hline
\multicolumn{1}{|l|}{$P_g^{\max}$}    &  Upper bound on power generation of $g$th generator/plant   \\ \hline
\multicolumn{2}{|l|}{\textit{C. Decision Variables}}       \\ \hline
\multicolumn{1}{|l|}{$p^G$}         & $g \times 1$ vector of baseline generations for meeting the demand $d$ with $p^R$        \\ \hline
\multicolumn{1}{|l|}{$p^R$}         & $u \times 1$ vector of the forecast outputs of renewable generators        \\ \hline
\multicolumn{1}{|l|}{$\lambda_i$}         & Participation factor of the $i$th power generator/plant on meeting the uncertain net load, where $0 \le \lambda_i \le 1$ and $\sum_{i \in \mathcal{R}}\lambda_i=1$  \\ \hline
\multicolumn{1}{|l|}{$\theta$}         & $n \times 1$ vector of phase angles of nodes/buses where, for the reference node/bus, $\theta_1=0$    \\ \hline
\multicolumn{2}{|l|}{\textit{D. Uncertain Variable}}       \\ \hline
\multicolumn{1}{|l|}{$\xi_d$}         &  Difference between the forecasted and real net demands of the $d$th load  \\ \hline
\end{tabular} \label{tab:nct}
\end{table}

\subsection{The PECO formulation for OPF under uncertainty}
When the uncertainty is considered, the real-time renewable generation $\hat{p}^R=p^R+\xi$, where $\xi \in \Upxi \subset \mathbb{R}^u$ is the uncertain component. The power generation $p^G$ needs to be adjusted in real-time such that the DCPF equation holds in real-time, i.e., $ B \hat{\theta} =A \hat{p}^G+C \hat{p}^R+d$ where the real-time phase angle $\hat{\theta}=\theta-\Delta \theta$. With a so-called affine control, i.e., $\hat{p}^G_i=p^G_i+\lambda_i\sum_{j \in \mathcal{R}}\xi_j$ ($\forall i \in \mathcal{G}$), the real-time DCPF equation is given as: 
   \vspace{-6pt}
\begin{equation} \label{rtDCPF}
   \vspace{-6pt}
    B (\theta-\Delta \theta) =A(p^G+e^{\rm T}\xi \lambda)+C(p^R+\xi)+d
\end{equation}
where $e$ is a $r \times 1$ identity vector. By comparing (\ref{rtDCPF}) to (\ref{DCPF}), we know that
   \vspace{-6pt}
\begin{equation}
   \vspace{-6pt}
    \Delta \theta=-\Breve{B}(Ae^{\rm T}\xi \lambda+C\xi) \nonumber
\end{equation}
where $\Breve{B}=
\begin{bmatrix}
0 & 0\\
0 & \hat{B}^{-1}
\end{bmatrix}$ and $\hat{B}$ is a $(n-1) \times (n-1)$ matrix obtained by removing the first row and column from $B$. Denote $\Delta \theta_i=-[\Breve{B}(Ae^{\rm T}\xi \lambda+C\xi)]_i$, the PECO formulation of OPF under uncertainty of renewable energy is given as
   \vspace{-6pt}
\begin{subequations} \label{PCCO-OPF} 
\begin{align}
\text{\textbf{p-OPF}:}\; \min_{p^G, \theta, \lambda} \quad & \mathbb{E}[c(p^G+e^{\rm T}y \lambda)]   \\
 \mathrm{s.t.} \quad  & (\text{\ref{DCPF}}) \\
 & \left|B_{ij}\left(\theta_i-[\Breve{B}(Ae^{\rm T}y \lambda+Cy)]_i-\theta_j+[\Breve{B}(Ae^{\rm T}y \lambda+Cy)]_j\right)\right|\\
 &\le P_{ij}^{\max},\; \forall\{i,j \} \in \mathcal{E} \nonumber \\
 &\sum_{i \in \mathcal{G}} \lambda_i=1 \\
 &   p_i^{\min} \le p^G_i+\lambda_ie^{\rm T}y \le p_i^{\max}, \forall i\in \mathcal{G} \\
 &  \forall y \in \{ y \in \Upxi \,|\,\mathbb{P}[\xi=y] \ge \alpha\},
\end{align} 
\end{subequations}
where $\mathbb{E}[c(p^G+e^{\rm T}y \lambda)]=\sum_{i \in \mathcal{G}}(c_{i,2}((p^G_i)^2+\mathbb{V}[e^{\rm T}y]\lambda_i^2)+c_{i,1}p^G_i+c_{i,0})=c(p^G)+\mathbb{V}[e^{\rm T}y]\sum_{i \in \mathcal{G}}(c_{i,2}\lambda_i^2)$, given that $\mathbb{E}[e^{\rm T}y]=0$ and $\mathbb{V}[\cdot]$ denotes variance.

Applying the findings in Sections 3 and 4, the DeDA($\mathcal{D}_\alpha^z$) and DeDA($\mathcal{D}_\alpha^\eta$) of (p-OPF) are given as follows, respectively.
   \vspace{-6pt}
\begin{subequations} \label{DDA1-OPF} 
\begin{align}
\text{\textbf{d-OPF}}&(\mathcal{D}^z_\alpha):\; \min_{p^G, \theta, \lambda} \;  c(p^G)+\left(\sum_{i \in \mathcal{G}}(c_{i,2}\lambda_i^2)\right)\frac{1}{z}\sum^{z}_{k=1}\left(e^{\rm T}\xi_{\rm d}^{(k)}\right)^2   \\
 \mathrm{s.t.} \;  & (\text{\ref{DCPF}}) \\
 & \left|B_{ij}\left(\theta_i-[\Breve{B}(Ae^{\rm T}\xi_{\rm d}^{(k)} \lambda+C\xi_{\rm d}^{(k)})]_i-\theta_j+[\Breve{B}(Ae^{\rm T}\xi_{\rm d}^{(k)} \lambda+C\xi_{\rm d}^{(k)})]_j\right)\right| \label{Flowlimit}\\
 &\le P_{ij}^{\max},\; \forall\{i,j \} \in \mathcal{E} \nonumber \\
 &\sum_{i \in \mathcal{G}} \lambda_i=1 \\
 &   p_i^{\min} \le p^G_i+\lambda_ie^{\rm T}\xi_{\rm d}^{(k)} \le p_i^{\max}, \forall i\in \mathcal{G} \label{Powerlimit} \\
 &  \forall \xi_{\rm d}^{(k)}  \in \mathcal{D}_\alpha^{z},
\end{align} 
\end{subequations}
\vspace{-28 pt}
\begin{subequations}
    \begin{align}
        \text{\textbf{d-OPF}}(\mathcal{D}^\eta_\alpha):\; \min_{p^G, \theta, \lambda} \;&  c(p^G)+\left(\sum_{i \in \mathcal{G}}(c_{i,2}\lambda_i^2)\right)\frac{1}{z}\sum^{z_\eta}_{k=1}\left(D^{\eta}_ke^{\rm T}\xi_{\rm d}^{(k)}\right)^2   \\
 \mathrm{s.t.} \;  & (\text{\ref{DCPF}}),\,\text{and}\, (\text{\ref{Flowlimit}})-(\text{\ref{Powerlimit}}) \\
 & \forall \xi_{\rm d}^{(k)}  \in \mathcal{D}_\alpha^{\eta}
    \end{align}
\end{subequations}
where $D^{\eta}_k$ is the number of data points in the $\eta$-vicinity of $\xi_{\rm d}^{(k)}$.

\subsection{Test Systems and Scenario Sets}
We used three representative test systems, i.e., IEEE 6, 39, and- 118-bus systems \cite{testcase}. The problem sizes of the p-OPFs for these test systems are provided in Table \ref{tab:testcase} (recall that $n$ and $r$ are the numbers of decision and uncertain variables respectively, and $m$ is the number of constraints that are contributed by one data point). In this numerical experiment, $\alpha=1\%$ is considered. Table \ref{tab:testcase} also tabulates the sizes of sets $\mathcal{D}$, $\mathcal{D}_\alpha$, $\mathcal{D}_\alpha^z$, and $\mathcal{D}_\alpha^\eta$ (recall that they are the original historical data set, the set of probable data, the embedded data sets of DeDA($\mathcal{D}_\alpha^z$) and DeDA($\mathcal{D}_\alpha^\eta$), respectively) for each case. Since both $x$ and $\xi$ are continuous, the above sets are determined following the path below: $$\mathcal{D} \xrightarrow[\text{to determine probable data}]{\text{Solve problem (\ref{jointprob2})}} \mathcal{D}_\alpha \xrightarrow[\text{with $z$ determined by (\ref{varrho})}]{\text{Randomly pick $z$ data points}} \mathcal{D}_\alpha^z \xrightarrow[\text{with $\eta=\zeta$}]{\text{ Algorithm 4.1}} \mathcal{D}_\alpha^\eta,$$ where we use the $\zeta$ determined in \ref{jointprob2} as the $\eta$ in Algorithm 4.1.

\begin{table}[h]
\footnotesize
\centering
\caption{p-OPF sizes of different test systems and the sizes of data sets that are used in these cases.}
\begin{tabular}{ccccccccccc}
\hline
\textbf{Case}     & \textit{n} & \textit{r} & \textit{m} & \textit{D} & $\zeta$ & $\alpha$  & $D_\alpha$ & $\rho$& $z$ & $z_{\eta=\zeta}$ \\ \hline
\textbf{IEEE-6}   & 9         & 2          & 16         & 1000      & 0.09 &0.05 & 685 &0.90  & 60       & 35  \\
\textbf{IEEE-39}   & 58          & 2          & 96         & 5000      & 0.12 & 0.01& 4459 &0.99 & 678       & 276  \\
\textbf{IEEE-118} & 155         & 10         & 315        & 10000      & 0.16 &0.01 & 9762 &0.99 & 773       & 368  \\ \hline
\end{tabular} \label{tab:testcase}
\end{table}


\subsection{Results and Analysis}
Recall that $x^*_{\rm D}(\mathcal{D}_\alpha)=x^*_{\rm P}$ under Assumption 1, we use DeDA($\mathcal{D}_\alpha$) as the reference for evaluating the performance of DeDA($\mathcal{D}_\alpha^z$) and DeDA($\mathcal{D}_\alpha^\eta$), i.e., the proposed DAs of PECO (\ref{PCCSP}). For each IEEE test case, the computational times and optimality gaps of DeDA($\mathcal{D}_\alpha^z$) and DeDA($\mathcal{D}_\alpha^\eta$) are compared in Table \ref{tab:results}. First, we can observed that the numerical results satisfy relation (4.1). Second, which is more important, DeDA($\mathcal{D}_\alpha^\eta$) is an accurate approximation to DeDA($\mathcal{D}_\alpha^z$) with a lower computational burden.

\begin{table}[h]
\footnotesize
\centering
\caption{Optimization results.}
\begin{tabular}{|c|c|c|c|c|c|}
\hline
\textbf{\begin{tabular}[c]{@{}c@{}}Test\\ system\end{tabular}} & \textbf{DeDA} & \textbf{\begin{tabular}[c]{@{}c@{}}Number of \\ constraints\end{tabular}} & \textbf{\begin{tabular}[c]{@{}c@{}}Computational\\ time (s)\end{tabular}} & \textbf{\begin{tabular}[c]{@{}c@{}}Objective\\ value (\$)\end{tabular}} & \textbf{\begin{tabular}[c]{@{}c@{}}Optimality\\ gap (\%)\end{tabular}} \\ \hline
\multirow{3}{*}{\textbf{IEEE-6}}                               & DeDA($\mathcal{D}_\alpha$)        & 10,960                                                                    & 0.66                                                                      & 2260.21                                                                 & --                                                                     \\ \cline{2-6} 
                                                               & DeDA($\mathcal{D}_\alpha^z$)        & 960                                                                       & 0.28                                                                      & 2258.73                                                                 & 0.0065\%                                                               \\ \cline{2-6} 
                                                               & DeDA($\mathcal{D}_\alpha^\eta$)        & 560                                                                       & 0.25                                                                      & 2257.71                                                                 & 0.1106\%                                                               \\ \hline
\multirow{3}{*}{\textbf{IEEE-39}}                              & DeDA($\mathcal{D}_\alpha$)            & 428,064                                                                   & 476.32                                                                    & 1876.67                                                                 & --                                                                     \\ \cline{2-6} 
                                                               & DeDA($\mathcal{D}_\alpha^z$)            & 65,088                                                                    & 161.81                                                                    & 1876.21                                                                 & 0.0025\%                                                               \\ \cline{2-6} 
                                                               & DeDA($\mathcal{D}_\alpha^\eta$)            & 26,496                                                                    & 72.54                                                                     & 1876.19                                                                 & 0.0026\%                                                               \\ \hline
\multirow{3}{*}{\textbf{IEEE-118}}                             & DeDA($\mathcal{D}_\alpha$)            & 3,075,030                                                                 & out of memory                                                             & --                                                                      & --                                                                     \\ \cline{2-6} 
                                                               & DeDA($\mathcal{D}_\alpha^z$)            & 243,465                                                                   & 4791.34                                                                   & 84901.03                                                                & --                                                                     \\ \cline{2-6} 
                                                               & DeDA($\mathcal{D}_\alpha^\eta$)            & 115,920                                                                   & 818.03                                                                    & 84838.03                                                                & 0.0074\%                                                               \\ \hline
\end{tabular} \label{tab:results}
\end{table}

\section{Conclusion}
For solving a specific PECO problem, we assume that a big data set $\mathcal{D}$ of historical measurements of the uncertain parameters $\xi$ is available. First, problem (\ref{jointprob1}) or (\ref{jointprob2}) is solve to determine the probable data points. With all non-probable data points removed, one can obtain set $\mathcal{D}_\alpha$. Second, $z$ data points are randomly selected from $\mathcal{D}_\alpha$ and stored in $\mathcal{D}_\alpha^z$, where $z$ is determined by a desired probability $\rho$ via $z=\underline{\varrho}^{-1}(\rho)$. Finally, $z_\eta$ data points are further selected from $\mathcal{D}_\alpha^z$ and stored in $\mathcal{D}_\alpha^{\eta}$ by an SDS algorithm if DeDA($\mathcal{D}_\alpha^z$) is still too large to compute. 

This research rethinks the entire procedure of solving optimization problems under uncertainty from logic modeling to solution algorithm design. The logic model (\ref{PCCSP}), i.e., PECO, defined in this paper is a novel alternative to the existing CCO. The PEC (\ref{PCC}) therein logically means that an optimal solution should be feasible to all probable realizations of the uncertain variables. Such a logical meaning grants the PECO a very high application value since it reflects the need of many engineering systems in terms of optimization under uncertainty. Since the existing solution methods are either inapplicable or inefficient to PECO, another key contribution of this paper lies in the novel solution paradigm which consists of data-embedded deterministication (as detailed in Section 3) and strategic data selection (as detailed in Section 4). With the proposed solution paradigm, PECO problems can be solved accurately with relatively low computational complexity. 

\appendix

\section{Proofs in Section 2}
\subsection{Proof of Proposition \ref{pro:feasisetpccsp}}
When $\alpha_1 \ge \alpha_2$, for an arbitrary outcome $\check{\xi}_{\rm s}$ of $\xi$, $\mathbb{P}[\xi=\check{\xi}_{\rm s}] \ge \alpha_2$ if $\mathbb{P}[\xi=\check{\xi}_{\rm s}] \ge \alpha_1$. In other words,  $\Omega_{\alpha_1} \subseteq \Omega_{\alpha_2}$. Then,
\begin{align}
    \mathcal{X}_{\rm P}(\alpha_2)&=\left\{x \in \mathbb{R}^n \middle\vert \begin{array}{l}
        g(x,y) \le 0,\;\forall y\in \Omega_{\alpha_1} \\
        g(x,y) \le 0,\;\forall y\in \Omega_{\alpha_2} \setminus \Omega_{\alpha_1}
    \end{array} \right\} \nonumber \\
    &=\{x \in \mathcal{X}_{\rm P}(\alpha_1) \mid g(x,y) \le 0,\;\forall y\in \Omega_{\alpha_2} \setminus \Omega_{\alpha_1} \} \nonumber
  \\
    & \subseteq \mathcal{X}_{\rm P}(\alpha_1). \nonumber
\end{align}

\subsection{Proof of Proposition \ref{pro:CD}} 
We consider the experiment that an outcome of $\xi$ is randomly extracted from the sample space $\Upxi$ and let $\check{\xi}_{\rm s}$ denote this arbitrary outcome. According to the probability theory, the probability of event $\mathbb{P}[\xi=\check{\xi}_{\rm s}]\ge \alpha$, i.e., $\mathbb{P}[\mathbb{P}[\xi=\check{\xi}_{\rm s}]\ge \alpha]$, is given as 
\begin{equation}
    \mathbb{P}[\mathbb{P}[\xi=\check{\xi}_{\rm s}]\ge \alpha]=\mathbb{P}[P(\check{\xi}_{\rm s})\ge \alpha]=\int\cdots\int_{\{\xi \in \Omega|P(\xi)\ge \alpha\}}P(\xi)d\xi_u\cdots d\xi_1,
\end{equation}
recalling that $\Omega$ is the mathematical expression of $\Upxi$. As a result, condition (\ref{CD}) indicates that $\mathbb{P}[P(\check{\xi}_{\rm s})\ge \alpha]=1-\beta$. Let $x_{\rm P} \in \mathcal{X}_{\rm P}$ denote an arbitrarily feasible solution of PECO, the logical meaning of PEC, i.e., $g(x_{\rm P},\check{\xi}_{\rm s}) \le 0$ holds if $\mathbb{P}[\xi=\check{\xi}_{\rm s}]\ge \alpha$, implies that the probability of event $g(x_{\rm P},\check{\xi}_{\rm s}) \le 0$ is not less than the probability of event $\mathbb{P}[\xi=\check{\xi}_{\rm s}]\ge \alpha$, i.e., $\mathbb{P}[g(x_{\rm P},\check{\xi}_{\rm s}) \le 0] \ge \mathbb{P}[P(\check{\xi}_{\rm s})\ge \alpha]=1-\beta$. Therefore, $x_{\rm P} \in \mathcal{X}_{\rm C}$ according to the logical meaning of the chance constraint, which implies that $\mathcal{X}_{\rm P} \subseteq \mathcal{X}_{\rm C}$.

The condition ``$P(\xi^{(a)}_{\rm s}) \ge P(\xi^{(b)}_{\rm s})$ when $\xi^{(a)}_{\rm s}$ and $\xi^{(b)}_{\rm s}$ are arbitrarily scenarios in $\mathcal{M}(x)$ and $\Omega \setminus \mathcal{M}(x)$ respectively" implies that there exists a probability $v$ such that $P(\xi^{(a)}_{\rm s}) \ge v$ and $P(\xi^{(b)}_{\rm s}) \le v$. Then, $\mathcal{M}(x)$ becomes
\begin{equation} \label{proofCD}
    \mathcal{M}^\prime(x)=\{\xi \in \Upxi\,|\,g(x,\xi) \le 0\;\text{and}\;P(\xi)\ge v\}.
\end{equation}
It's not hard to know that, in (\ref{proofCD}), $v = \alpha$ under condition (\ref{CD}). Let $\hat{\xi}_{\rm s}$ be an arbitrary realization of $\xi$ in $\mathcal{M}^\prime(x)$, we have $\mathbb{P}[\xi=\hat{\xi}_{\rm s}]\ge \alpha$. Let $x_{\rm C}$ denote an arbitrarily feasible solution of CCO that satisfies the above condition, we have $g(x_{\rm C},\hat{\xi}_{\rm s}) \le 0$, which implies that $x_{\rm C} \in \mathcal{X}_{\rm P}$ and, namely, $\mathcal{X}_{\rm P} \supseteq \mathcal{X}_{\rm C}$. Therefore, we have $\mathcal{X}_{\rm P} = \mathcal{X}_{\rm C}$.

\subsection{Proof of Proposition \ref{RDA}} 
Considering the experiment of randomly generating $N$ i.i.d. samples, it is of high probability that these samples belong to scenarios whose probabilities are not less than $1/N$. When $N \le \frac{1}{\alpha}$, it's of high probability that these samples belong to scenarios whose probabilities are not less than $\alpha$. Note that $\mathcal{X}_{\rm S}$ is a set in $x$-space specified by $N$ samples whose probabilities are not less than $\alpha$, we have  $\mathcal{X}_{\rm P} \subseteq \mathcal{X}_{\rm S}$ since $\mathcal{X}_{\rm P}$ is specified by all possible scenarios whose probabilities are not less than $\alpha$. Given that $N \le \frac{1}{\alpha}$ is a relatively small number (for example $N=100$ when $\alpha=0.01$), it's of high probability that $\mathcal{X}_{\rm P} \subset \mathcal{X}_{\rm S}$. When $N > \frac{1}{\alpha}$, a scenario whose probability is less than $\alpha$ may be included in the $N$ samples, which results in that an element in $\mathcal{X}_{\rm P}$ may not be feasible to $\mathcal{X}_{\rm S}$. Hence, $\mathcal{X}_{\rm P} \neq \mathcal{X}_{\rm S}$ is high-probability event under this condition.

\section{Proofs in Section 3}
\subsection{Proof of Proposition \ref{pro:underlying}}
 If we number the data points in $\mathcal{D}_1$ as $\{\xi_{\rm d}^{(i)}\, (i=1,\ldots,I)\}$, we can renumber the data points in $\mathcal{D}_2$ as $\{\xi^{(i,j)}_{\rm d}\, (i =1,\ldots,I;\,j=1,\ldots,J_i)\}$. Note that $\xi^{(i,1)}_{\rm d}=,\ldots,=\xi^{(i,J_i)}_{\rm d}$ ($\forall i =1,\ldots,I$) according to the properties of a multiset and its underlying set. Then, DeDA($\mathcal{D}_2$) can be reformulated as
 \begin{subequations} \label{itm1}
\begin{align}
 \min_{x} \quad & f(x)  \\
 \mathrm{s.t.} \quad  &g(x,\xi^{(i,1)}_{\rm d}) \le 0 \label{constr1_itm1} \\
 &g(x,\xi^{(i,j)}_{\rm d}) \le 0, \label{constr2_itm1} 
\end{align} 
\end{subequations}
where $i =1,\ldots,I$ and $j=1,\ldots,J_i$. Since $\xi^{(i,1)}_{\rm d}=,\ldots,=\xi^{(i,z_i)}_{\rm d}$ ($\forall i =1,\ldots,I$), constraints in (\ref{constr2_itm1}) are redundant and can be removed without affecting the solution. Without constraints (\ref{constr2_itm1}), problem (\ref{itm1}) is exactly the DeDA($\mathcal{D}_1$). Namely, $\mathcal{X}_{\rm D}(\mathcal{D}_1)=\mathcal{X}_{\rm D}(\mathcal{D}_2)$.

\subsection{Proof of Proposition \ref{pro:adpfs}}

Since $\mathcal{S}_\alpha^{\forall}$ is defined as the set that contains all possible scenarios of $\Upxi_\alpha$, we have $\mathfrak{U}[\mathcal{D}_\alpha^z] \subseteq \mathcal{S}_\alpha^{\forall}$ and, consequently $\mathcal{X}_{\rm D}(\mathcal{D}_\alpha^z)=\mathcal{X}_{\rm D}(\mathfrak{U}[\mathcal{D}_\alpha^z]) \supseteq \mathcal{X}_{\rm D}(\mathcal{S}_\alpha^{\forall})$ according to Propositions \ref{pro:underlying} and \ref{pro:feasisetpccsp} (Note that $\mathcal{S}_\alpha^{\forall}$ is not a multiset). According to Definition \ref{def:equivalent}, we have $\mathcal{X}_{\rm D}(\mathcal{S}_\alpha^{\forall})=\mathcal{X}_{\rm P}$ and, consequently, $\mathcal{X}_{\rm D}(\mathcal{D}_\alpha^z) \supseteq \mathcal{X}_{\rm P}$. Condition $\mathcal{B}_{\forall,\alpha}^{\rm FS} \subseteq \mathcal{D}_\alpha^z$ indicates that $\mathcal{X}_{\rm D}(\mathcal{B}_{\forall,\alpha}^{\rm FS}) \supseteq \mathcal{X}_{\rm D}(\mathcal{D}_\alpha^z)$. According to Definition \ref{def:adpfs}, we have $\mathcal{X}_{\rm D}(\mathcal{B}_{\forall,\alpha}^{\rm FS})=\mathcal{X}_{\rm D}(\mathcal{S}_\alpha^{\forall})=\mathcal{X}_{\rm P}$, which means $\mathcal{X}_{\rm P}\supseteq \mathcal{X}_{\rm D}(\mathcal{D}_\alpha^z)$. Then, we have $\mathcal{X}_{\rm D}(\mathcal{D}_\alpha^z) = \mathcal{X}_{\rm P}$ if $\mathcal{B}_{\forall,\alpha}^{\rm FS} \subseteq \mathcal{D}_\alpha^z$.

If $\mathcal{X}_{\rm D}(\mathcal{D}_\alpha^z) = \mathcal{X}_{\rm P}$, we have $\mathcal{X}_{\rm D}(\mathcal{D}_\alpha^z) =\mathcal{X}_{\rm D}(\mathcal{S}_\alpha^{\forall})$. Assuming $\mathcal{B}_{\forall,\alpha}^{\rm FS} \not \subseteq \mathcal{D}_\alpha^z$ and letting $\Tilde{\xi}_{\rm d} \not \in \mathcal{D}_\alpha^z$ be an element in $\mathcal{B}_{\forall,\alpha}^{\rm FS}$ and $\mathcal{D}_\alpha^\prime=\mathcal{D}_\alpha^z \bigcup \Tilde{\xi}_{\rm d}$, there are two possible relations between $\mathcal{X}_{\rm D}(\mathcal{D}_\alpha^z)$ and $\mathcal{X}_{\rm D}(\mathcal{D}_\alpha^\prime)$, i.e., $\mathcal{X}_{\rm D}(\mathcal{D}_\alpha^z)=\mathcal{X}_{\rm D}(\mathcal{D}_\alpha^\prime)$ and $\mathcal{X}_{\rm D}(\mathcal{D}_\alpha^z)\supset\mathcal{X}_{\rm D}(\mathcal{D}_\alpha^\prime)$, respectively. If $\mathcal{X}_{\rm D}(\mathcal{D}_\alpha^z)=\mathcal{X}_{\rm D}(\mathcal{D}_\alpha^\prime)$, we have $\mathcal{X}_{\rm D}(\mathcal{D}_\alpha^z)=\mathcal{X}_{\rm D}(\mathcal{D}_\alpha^\prime)=\mathcal{X}_{\rm D}(\mathcal{S}_\alpha^{\forall})$, which implies that removing $\Tilde{\xi}_{\rm d}$ from $\mathcal{B}_{\forall,\alpha}^{\rm FS}$ (or $\mathcal{S}_\alpha^{\forall}$) does not affect the feasible set. This feature contradicts Definition \ref{def:adpfs} that $\mathcal{B}_{\forall,\alpha}^{\rm FS}$ is the smallest subset of $\mathcal{S}_\alpha^{\forall}$ which satisfies condition (\ref{condition_adpfs}). Moreover, we know that $\mathcal{X}_{\rm D}(\mathcal{D}_\alpha^\prime) \supseteq \mathcal{X}_{\rm D}(\mathcal{S}_\alpha^{\forall})$ since $\mathcal{D}_\alpha^\prime \subseteq \mathcal{S}_\alpha^{\forall}$. If $\mathcal{X}_{\rm D}(\mathcal{D}_\alpha^z)\supset\mathcal{X}_{\rm D}(\mathcal{D}_\alpha^\prime)$, it contradicts the precondition of $\mathcal{X}_{\rm D}(\mathcal{D}_\alpha^z) = \mathcal{X}_{\rm P}$. Hence, we have $\mathcal{B}_{\forall,\alpha}^{\rm FS} \subseteq \mathcal{D}_\alpha^z$ if $\mathcal{X}_{\rm D}(\mathcal{D}_\alpha^z) = \mathcal{X}_{\rm P}$.

\subsection{Proof of Proposition \ref{pro:adpos2}}
Recall that $\mathcal{X}^*_{\rm D}(\Box)$ and $x^*_{\rm D}(\Box)$ are the feasible set and optimal solution, respectively, of DeDA($\Box$), and $\mathcal{B}_{\Box}^{\rm FS}$ is the set of boundary-forming data points of $\mathcal{X}^*_{\rm D}(\Box)$. Definition \ref{def:adpfs} implies $x^*_{\rm D}(\mathcal{B}_{\Box}^{\rm FS})=x^*_{\rm D}(\Box)$ since $\mathcal{X}^*_{\rm D}(\mathcal{B}_{\Box}^{\rm FS})$ $=\mathcal{X}^*_{\rm D}(\Box)$. Definition \ref{def:adpos} means that $\mathcal{B}_{\Box}^{\rm OS}$ contains all boundary-forming data points of $\mathcal{X}^*_{\rm D}(\Box)$ which are active to $x^*_{\rm D}(\Box)$. In other words, the data points in the complementary of $\mathcal{B}_{\Box}^{\rm OS}$ with respect to $\mathcal{B}_{\Box}^{\rm FS}$ are inactive to $x^*_{\rm D}(\Box)$. Then, it is not hard to know that, removing all data points in $\mathcal{B}_{\Box}^{\rm FS} \setminus \mathcal{B}_{\Box}^{\rm OS}$, the optimization problem reduce from DeDA($\mathcal{B}_{\Box}^{\rm FS}$) to DeDA($\mathcal{B}_{\Box}^{\rm OS}$) while the optimal solution does not change. Hence, $x^*_{\rm D}(\mathcal{B}_{\Box}^{\rm OS})=x^*_{\rm D}(\mathcal{B}_{\Box}^{\rm FS})=x^*_{\rm D}(\Box)$. Recalling that $\mathcal{B}_{\forall,\alpha}^{\rm OS}$ denotes the set of boundary-forming data points at $x^*_{\rm D}(\mathcal{S}_\alpha^{\forall})$, we have $x^*_{\rm D}(\mathcal{B}_{\forall,\alpha}^{\rm OS})=x^*_{\rm D}(\mathcal{S}_\alpha^{\forall})$. If $\mathcal{B}_{\forall,\alpha}^{\rm OS} \subseteq \mathcal{D}_\alpha^z$, i.e.,  $\mathcal{D}_\alpha^z$ contains all boundary-forming data points of $x^*_{\rm D}(\mathcal{S}_\alpha^{\forall})$, we have $\mathcal{B}_{\forall,\alpha}^{\rm OS} \subseteq \mathfrak{U}[\mathcal{D}_\alpha^z]$ since $\mathcal{B}_{\forall,\alpha}^{\rm OS}$ is not a multiset. As a result, we have $(\mathcal{S}_\alpha^{\forall}\setminus \mathfrak{U}[\mathcal{D}_\alpha^z])\subseteq(\mathcal{S}_\alpha^{\forall} \setminus \mathcal{B}_{\forall,\alpha}^{\rm OS})$. With all data points in $(\mathcal{S}_\alpha^{\forall}\setminus \mathfrak{U}[\mathcal{D}_\alpha^z])$ removed, the optimization problem DeDA($\mathcal{S}_\alpha^{\forall}$) reduces to DeDA($\mathfrak{U}[\mathcal{D}_\alpha^z]$). It suffices to know that the optimal solution $x^*_{\rm D}(\mathcal{S}_\alpha^{\forall})$ does not change during this process since all data points in $\mathcal{S}_\alpha^{\forall} \setminus \mathcal{B}_{\forall,\alpha}^{\rm OS}$ are inactive to $x^*_{\rm D}(\mathcal{S}_\alpha^{\forall})$. In other words, we have $x^*_{\rm D}(\mathfrak{U}[\mathcal{D}_\alpha^z])=x^*_{\rm D}(\mathcal{S}_\alpha^{\forall})$ and consequently, $x^*_{\rm D}(\mathcal{D}_\alpha^z)=x^*_{\rm P}$ since $x^*_{\rm D}(\mathcal{D}_\alpha^z)=x^*_{\rm D}(\mathfrak{U}[\mathcal{D}_\alpha^z])$ and $x^*_{\rm P}=x^*_{\rm D}(\mathcal{S}_\alpha^{\forall})$.


\subsection{Proof of Theorem \ref{thm:varrho}}
Let $\xi_{\rm d,B}^{(i)}$ ($i=1,\ldots,B_{\forall,\alpha}^{\rm OS}$) denote the boundary-forming data points for $x^*_{\rm D}(\mathcal{S}_\alpha^{\forall})$ and $N_i$ be the number of $\xi_{\rm d,B}^{(i)}$ in $\mathcal{D}_\alpha$, we consider the following events:
\begin{itemize}[leftmargin=*]
    \item E$^\circ$: when a data point is randomly selected from $\mathcal{D}_\alpha$, it is $\xi_{\rm d,B}^{(i)}$, i.e., one of the $B_{\forall,\alpha}^{\rm OS}$ boundary-forming data points;
    \item E$^*$: when $z$ data points are randomly selected from $\mathcal{D}_\alpha$, at least one of each of the $B_{\forall,\alpha}^{\rm OS}$ active data points is selected;
    \item E$_i$ ($ i=1,\ldots,B_{\forall,\alpha}^{\rm OS}$): when $z$ data points are randomly selected from $\mathcal{D}_\alpha$, no $\xi_{\rm d,B}^{(i)}$ is selected.
\end{itemize}
 It suffices to show that $\overline{\text{E}^*}=\bigcup_{i=1}^{B_{\forall,\alpha}^{\rm OS}}\text{E}_i$ which is the complement of $\text{E}^*$. Then, we have $\mathbb{P}[\text{E}^*]=1-\mathbb{P}[\overline{\text{E}^*}]=1-\mathbb{P}[\bigcup_{i=1}^{B_{\forall,\alpha}^{\rm OS}}\text{E}_i]$ and 

\begin{align}
\mathbb{P}\left[\bigcup_{i=1}^{B_{\forall,\alpha}^{\rm OS}}\text{E}_i\right]=&\sum_{i=1}^{B_{\forall,\alpha}^{\rm OS}}\mathbb{P}[\text{E}_i]-\sum_{i=1}^{B_{\forall,\alpha}^{\rm OS}}\sum_{j>i}^{B_{\forall,\alpha}^{\rm OS}}\mathbb{P}[\text{E}_i\cap \text{E}_j]+\sum_{i=1}^{B_{\forall,\alpha}^{\rm OS}}\sum_{j>i}^{B_{\forall,\alpha}^{\rm OS}}\sum_{k>j}^{B_{\forall,\alpha}^{\rm OS}}\mathbb{P}[\text{E}_i\cap \text{E}_j\cap \text{E}_k]\\
    &-\cdots+(-1)^{(B_{\forall,\alpha}^{\rm OS}-1)}\mathbb{P}\left[\bigcap_{i=1}^{B_{\forall,\alpha}^{\rm OS}}\text{E}_i\right]. \nonumber
\end{align}

It's not hard to know that $\mathbb{P}[\text{E}_i]$, $\mathbb{P}[\text{E}_i\cap \text{E}_j]$, $\mathbb{P}[\text{E}_i\cap \text{E}_j\cap \text{E}_k]$, $\cdots$, and $\mathbb{P}[\bigcap_{i=1}^{B_{\forall,\alpha}^{\rm OS}}\text{E}_i]$ follow the hypergeometric distribution \cite{rice2006mathematical}. Hence, we have
\begin{equation}
    \mathbb{P}[K]= \mathbb{P}\left[\bigcap_{i=1}^K\text{E}_i\right]= \frac{\binom{D_\alpha-\sum_{i=1}^KN_i}{z}}{\binom{D_\alpha}{z}},\; K=1,\ldots,B_{\forall,\alpha}^{\rm OS}
\end{equation}
and
\begin{equation}
    \mathbb{P}[\text{E}^*]=1-\mathbb{P}\left[\bigcup_{i=1}^{B_{\forall,\alpha}^{\rm OS}}\text{E}_i\right]=1-\sum_{K=1}^{B_{\forall,\alpha}^{\rm OS}}\left((-1)^{K-1}\binom{B_{\forall,\alpha}^{\rm OS}}{K}\mathbb{P}[K]\right).
\end{equation}
Given that all the boundary-forming data points should belong to $\mathcal{S}_\alpha^{\forall}$, which implies $\mathbb{P}[\xi=\xi_{\rm d,B}^{(i)}] \ge \alpha$ ($ i=1,\ldots,B_{\forall,\alpha}^{\rm OS}$), the expectation of the times that $\xi_{\rm d,B}^{(i)}$ appears in $\mathcal{D}$ is not less than $\alpha D$, i.e., $\mathbb{E}[N_i]\ge \alpha D$. When $D$ is big, the Central Limit Theorem indicates that $N_i\approx \mathbb{E}[N_i]$. Hence, we consider $N_i\ge \alpha D$ ($ i=1,\ldots,B_{\forall,\alpha}^{\rm OS}$). Moreover, since $B_{\forall,\alpha}^{\rm OS} \le \bar B_{\forall,\alpha}^{\rm OS}$, we have 
\begin{equation} \label{thmvarrho}
    \mathbb{P}[\text{E}^*] \ge 1-\sum_{K=1}^{\bar B_{\forall,\alpha}^{\rm OS}}\left((-1)^{K-1}\binom{\bar B_{\forall,\alpha}^{\rm OS}}{K}\frac{\binom{D_\alpha-K\alpha D}{z}}{\binom{D_\alpha}{z}}\right)=\underline{\varrho}.
\end{equation}
According to Proposition \ref{pro:adpos2}, event E$^*$ implies that $x^*_{\rm D}(\mathcal{D}_\alpha^z)=x^*_{\rm D}(\mathcal{S}_\alpha^{\forall})$ and, consequently $x^*_{\rm D}(\mathcal{D}_\alpha^z)=x^*_{\rm P}$. Therefore, $\varrho =\mathbb{P}[x^*_{\rm D}(\mathcal{D}_\alpha^z)=x^*_{\rm P}]=\mathbb{P}[\text{E}^*]$ and the $\underline{\varrho}$ in (\ref{thmvarrho}) is a lower bound of $\varrho$.

\subsection{Proof of Proposition \ref{pro:adpos1}}
To prove this proposition, we utilize the finding on the support constraint, as defined in \cite{calafiore2005uncertain}, whose removal changes the solution of the optimization problem. According to Proposition 1 in \cite{calafiore2005uncertain}, for a convexly constrained linear program, the number of support constraints is at most $n$ (read \cite{calafiore2005uncertain} for a detailed proof). Let $x^*_{\rm D}(\Box)$ denote the optimal solution of DeDA (\ref{DDA}), by comparing Definition \ref{def:adpos} of this paper and the definition of the support constraint, we know that the number of boundary-forming constraints is less than or equal to that of the support constraints. Given that a boundary-forming data point of $x^*_{\rm D}(\Box)$ contributes at least one boundary-forming constraint to $x^*_{\rm D}(\Box)$ (note that a data point contributes $m$ constraints to the DeDA). As such, we have ``$n \ge$ the number of support constraints $\ge$ the number of boundary-forming constraints $\ge$ the number of boundary-forming data points." 

\section{Proofs in Section 4}
\subsection{Proof of Theorem \ref{pro:continuous}}
When $\xi \in \Upxi \subset \mathbb{R}^u$, it's straightforward to know that $\mathcal{D}_\alpha^\eta=\mathcal{D}_\alpha^z$ if $\eta=0$. Hence, $z_\eta=z$ when $\eta=0$. When $\eta \ge \bar \eta$, only one data point (i.e., $\xi^{(0)}_{\rm d}$) is selected from $\mathcal{D}_\alpha^z$ and stored in $\mathcal{D}_\alpha^\eta$ since no other data points satisfy condition (\ref{agr2}). In the $u$-dimensional Euclidean space, the $\eta$-vicinity of a data point is actually an $u$-ball whose volume is $V_u^\eta=2(2 \pi)^{\frac{u-1}{2}} \eta^u/u!!$ \cite{lozier2003nist}, where $u!!$ is the double factorial of $u$. Let $\mathcal{A} \subset \Upxi \subset \mathbb{R}^u$ be the smallest continuous set (i.e., convex hull) that contains all the data points in $\mathcal{D}^z_\alpha$ and $V_u^{\mathcal{A}}$ denote the $u$-dimensional Euclidean volume of $\mathcal{A}$, we divide $V_u^{\mathcal{A}}$ by $V_u^\eta$ and denote it as $\bar z_\eta$:
 \begin{equation} \label{Neta1}
    \bar z_\eta =\frac{V_u^{\mathcal{A}}}{V_u^\eta}= \frac{u!! V_u^{\mathcal{A}} }{2\eta^u(2 \pi)^{\frac{u-1}{2}}},
\end{equation}
 where $0< \eta < \bar \eta$. One can consider Algorithm \ref{alg:sds} as ``packing $z_\eta$ $u$-balls in the $u$-polytope $\mathcal{A}$." It's straightforward to know that $z_\eta<\bar z_\eta$ due to the existence of ``gaps" among the $z_\eta$ $u$-balls. In other words, we can consider $\bar z_\eta$ an upper estimate of $z_\eta$ for $0< \eta < \bar \eta$, and relation (\ref{Neta1}) implies $\bar z_\eta \propto 1/\eta^u$. 

 Since the positions of the $u$-balls are randomly chosen when Algorithm \ref{alg:sds} packs these $u$-balls in $\mathcal{A}$, the needed number $z_\eta$ of $u$-balls for filling out $\mathcal{A}$ would slightly vary even when the radius $\eta$ of these $u$-balls is fixed. Therefore, $z_\eta$ is a random variable which is parameterized by $\eta$. It is also straightforward to know that, if the radius $\eta$ of these $u$-balls is bigger, the number $z_\eta$ of $u$-balls that can be packed in $\mathcal{A}$ is most likely less. Now, let's consider two sets of random experiments where Algorithm \ref{alg:sds} fills $\mathcal{A}$ with $u$-balls of radius $\eta_{(1)}$ in the first set of experiments and with $u$-balls of radius $\eta_{(2)}$ in the second. Further note that $\eta_{(1)} < \eta_{(2)}$ and, in each experiment, the positions of $r$-balls are randomly selected. It's not difficult to know that the expectation of $z_{\eta_{(1)}}$ should not be less than that of $z_{\eta_{(2)}}$, i.e., $\mathbb{E}[z_{\eta_{(1)}}] \ge \mathbb{E}[z_{\eta_{(2)}}]$. Then,
 \begin{equation}
     \frac{\mathbb{E}[z_{\eta_{(1)}}] - \mathbb{E}[z_{\eta_{(2)}}]}{\eta_{(1)} - \eta_{(2)}} \le 0. \nonumber
 \end{equation}
 The limit of the left-hand-side of the above inequality as $|\eta_{(1)} - \eta_{(2)}| \rightarrow 0$ yields $\partial (\mathbb{E}[z_\eta])/\partial \eta \le 0$.

\subsection{Proof of Theorem \ref{thm:rho2}}
 Recall that the number of boundary-forming constraints at $x^*_{\rm D}(\mathcal{D}_\alpha^z)$ was assumed to be $B_c$ and are denoted as $\Tilde{g}(x,\xi_{\rm d,B}) \le 0$ ($\Tilde{g}:\mathbb{R}^{n+(u\times B_c)} \rightarrow \mathbb{R}^{B_c}$), we have $\Tilde{g}(x^*_{\rm D}(\mathcal{D}_\alpha^z),\xi_{\rm d,B})=0$ while condition 1 asserts that $\Tilde{g}(x^*_{\rm D}(\mathcal{D}_\alpha^{\eta}),\acute{\xi}_{\rm d}) = 0$. Proposition \ref{pro:ift} indicates the existence of a vector-valued function $x=h(\Tilde{\xi})$ which is equivalent to $\Tilde{g}(x,\Tilde{\xi})=0$ in the vicinity of ($x^*_{\rm D}(\mathcal{D}_\alpha^z),\xi_{\rm d,B}$) under the conditions therein. Let $\Delta x=x^*_{\rm D}(\mathcal{D}_\alpha^{\eta})-x^*_{\rm D}(\mathcal{D}_\alpha^z)$ and $\Delta \Tilde{\xi}=\acute{\xi}_{\rm d}-\xi_{\rm d,B}$, the Taylor series of $x=h(\Tilde{\xi})$ at $\xi_{\rm d,B}$ is 
\begin{equation}
    \Delta x=H(\xi_{\rm d,B})\Delta \Tilde{\xi}+\frac{1}{2!}\Delta \Tilde{\xi}^{\rm T}H^\prime(\xi_{\rm d,B})\Delta \Tilde{\xi} + \frac{1}{3!}\Delta \Tilde{\xi}^{\rm T}H^{\prime\prime}(\xi_{\rm d,B})(\Delta \Tilde{\xi})^2+\cdots. \label{taylor}
\end{equation}
Then, we have
\begin{equation} \label{taylor1}
     \|\Delta x\| \le\|H(\xi_{\rm d,B})\| \|\Delta \Tilde{\xi}\|+\frac{1}{2!}\|H^\prime(\xi_{\rm d,B})\| \|\Delta \Tilde{\xi}\|^2 + \frac{1}{3!}\|H^{\prime\prime}(\xi_{\rm d,B})\| \|\Delta \Tilde{\xi}\|^3+\cdots
\end{equation}
Note that $B_c=n$ under the conditions in Proposition \ref{pro:ift} and $\Delta \Tilde{\xi}=[\Delta \Tilde{\xi}_1^{\rm T},\ldots, \Delta \Tilde{\xi}_{n}^{\rm T}]^{\rm T}$, where $\|\Delta \Tilde{\xi}_i\| \le 2\eta$ ($i=1,\ldots,n$) according to Algorithm \ref{alg:sds}, we have $\|\Delta \Tilde{\xi}\| \le 2\sqrt{n} \eta$. Therefore, from (\ref{taylor1}), we have:
\begin{equation} \label{deltax}
    \|\Delta x\| \le 2\sqrt{n}\|H(\xi_{\rm d,B})\|  \eta+\frac{2^2n}{2!}\|H^\prime(\xi_{\rm d,B})\|\eta^2+\frac{2^3\sqrt{n^3}}{3!}\|H^{\prime\prime}(\xi_{\rm d,B})\|\eta^3+\cdots
\end{equation}
Then, we further have:
\begin{align}
    \varphi&=1-\frac{\|\Delta x\|}{\| x^*_{\rm D}(\mathcal{D}^z_{\alpha})\|} \nonumber \\
    &\ge1-\frac{2\sqrt{n}\|H(\xi_{\rm d,B})\|\eta+\frac{2^2n}{2!}\|H^\prime(\xi_{\rm d,B})\|\eta^2+\frac{2^3\sqrt{n^3}}{3!}\|H^{\prime\prime}(\xi_{\rm d,B})\|\eta^3+\cdots}{\| x^*_{\rm D}(\mathcal{D}^z_{\alpha})\|}\\
    &=\underline{\varphi}(\eta). \nonumber
\end{align}

When the boundary-forming data points of $x^*_{\rm D}(\mathcal{D}_\alpha^{\eta})$ are not exactly $\acute{\xi}_{\rm d}^{(i)}$ ($\forall i=1,\dots,B_c$), we let $\acute{x}$ denote the solution of $\Tilde{g}(x,\acute{\xi}_{\rm d}) = 0$. The condition, that $g$ is convex in $x$, implies the convexity of $\mathcal{X}_{\rm D}$'s. Due to the convexity of $f$, it's not hard to know that $x^*_{\rm D}(\mathcal{D}_\alpha^{\eta})=x^*_{\rm D}(\mathcal{D}_\alpha^z)$ if $x^*_{\rm D}(\mathcal{D}_\alpha^{\eta}) \in \mathcal{X}^*_{\rm D}(\mathcal{D}_\alpha^z)$. Now, let's consider the situation of $x^*_{\rm D}(\mathcal{D}_\alpha^{\eta}) \notin \mathcal{X}^*_{\rm D}(\mathcal{D}_\alpha^z)$. Let $\acute{\mathcal{D}}=\{\acute{\xi}_{\rm d}^{(i)}\;( i=1,\dots,B_c)\}$, it's not hard to know that $\acute{x}$ is the optimal solution of DeDA($\acute{\mathcal{D}}$) and $\acute{\mathcal{D}} \subseteq \mathcal{D}_\alpha^{\eta}$. Then, we have $\mathcal{X}^*_{\rm D}(\acute{\mathcal{D}}) \supseteq \mathcal{X}^*_{\rm D}(\mathcal{D}_\alpha^{\eta})\supseteq \mathcal{X}^*_{\rm D}(\mathcal{D}_\alpha^z)$, $\acute{x} \notin \mathcal{X}^*_{\rm D}(\mathcal{D}_\alpha^z)$, and $f(\acute{x})\le f(x^*_{\rm D}(\mathcal{D}_\alpha^{\eta}))$. If we let $d_1=\Delta x_{(1)}=x^*_{\rm D}(\mathcal{D}_\alpha^{\eta})-x^*_{\rm D}(\mathcal{D}_\alpha^z)$ and $d_2=\Delta x_{(2)}=\acute{x}-x^*_{\rm D}(\mathcal{D}_\alpha^z)$, both $d_1$ and $d_2$ satisfy condition (\ref{outward}) since a direction that satisfies (\ref{outward}) is an ``outward" pointing direction. The condition $f(\acute{x})\le f(x^*_{\rm D}(\mathcal{D}_\alpha^{\eta})) \le f(x^*_{\rm D}(\mathcal{D}_\alpha^z))$ implies $\|\Delta x_{(1)}\| \le \| \Delta x_{(2)}\|$ since $f$ is radially non-increasing in the vicinity of $x^*_{\rm D}(\mathcal{D}_\alpha^z)$ alone any direction $d$ that satisfies (\ref{outward}). In the previous paragraph, we already showed that $\Delta x_{(2)}$ satisfies condition (\ref{deltax}). Hence, the lower bound of $\varphi$ given in (\ref{thm4.3}) is also valid under condition 2.

 \subsection{Proof of Corollary \ref{cor:Neta}}
 Recall that all elements in a subset $\mathcal{D}^z_{\alpha,i}$ ($i \in \mathcal{I}$) of $\mathcal{D}_\alpha^z$ have the same integer part. With the integer part being fixed, steps 2-7 in Algorithm \ref{alg:mix} on $\mathcal{D}^z_{\alpha,i}$ is equivalent to Algorithm \ref{alg:sds} on $\mathcal{D}^z_{\alpha}$. Therefore, one can apply the assertions in Theorem \ref{pro:continuous} to each $i \in \mathcal{I}$. Then, we have  
   \begin{equation}  
    z_{\eta,i} \begin{cases} =z_i &  \eta_i =0 \\\le \Bar{z}_{\eta,i} & 0< \eta_i < \bar \eta_i \\ =1 & \eta_i \ge \bar \eta_i \end{cases},\;\Bar{z}_{\eta,i} \propto \frac{1}{\eta_i^{r_2}},\;\text{and}\;     \frac{\partial (\mathbb{E}[z_{\eta,i}])}{\partial \eta_i}   \le 0, \; i \in \mathcal{I}. \nonumber
    \end{equation}
Since $z=\sum_{i\in \mathcal{I}}z_i$ and $z_\eta=\sum_{i\in  \mathcal{I}}z_{\eta,i}$, we have (\ref{cor_con_1}) and (\ref{cor_con_2}).

\subsection{Proof of Corollary \ref{cor:rho2}}
 When $\eta_i$ ($\forall i \in \mathcal{I}$) are small, $\xi_{\rm d,B}$, $\xi_{\rm d,B^\prime}$, and $\acute{\xi}_{\rm d}$ share the same integer parts. As a result, the situation reduces to that of continuous $\xi$ and the assertions in Theorem \ref{thm:rho2} are applicable to this situation. We don't need to consider the subsets $\mathcal{D}^z_{\alpha,i}$ ($i \in \mathcal{I}$) of $\mathcal{D}_\alpha^z$ which do not contain active data points of $x^*_{\rm D}(z)$ since removing any inactive data point will not impact the optimal solution. Under Algorithm \ref{alg:mix}, we have $\|\Delta \xi_i\| \le 2\eta_i$ ($i=1,\ldots,B_c$), which implies that $\|\Delta \xi\| \le 2\sqrt{\eta_1^2+\ldots+\eta_{B_c}^2}= \hat{\eta}$. Then, applying the conclusions in Theorem \ref{thm:rho2}, we have (\ref{cor_rho2}). 

\subsection{Proof of Lemma \ref{pro:discrete}}
When $\xi \in \Upxi \subset \mathbb{Z}^u$, we already know that $\mathcal{D}_\alpha^{\eta}=\mathfrak{U}[\mathcal{D}_\alpha^z]$. We also know that $\mathcal{X}_{\rm D}(\mathfrak{U}[\Box])=\mathcal{X}_{\rm D}(\Box)$ from Proposition \ref{pro:underlying}. Therefore, we have $\mathcal{X}_{\rm D}(\mathcal{D}_\alpha^{\eta})=\mathcal{X}_{\rm D}(\mathcal{D}_\alpha^z)$ and $x^*_{\rm D}(\mathcal{D}_\alpha^{\eta})=x^*_{\rm D}(\mathcal{D}_\alpha^z)$. Recall that $x^*_{\rm D}(\mathcal{D}_\alpha)$ denotes the optimal solution of DeDA($\mathcal{D}_\alpha$), Assumption 1 implies that $x^*_{\rm D}(\mathcal{D}_\alpha)=x^*_{\rm P}$. If Algorithm \ref{alg:mix} uses $\mathcal{D}_\alpha$ as input data set, it's not hard to know that the resulting $\mathcal{D}_\alpha^{\eta}=\mathfrak{U}[\mathcal{D}_\alpha]$. Following the analysis above, it suffices to have $x^*_{\rm D}(\mathcal{D}_\alpha^{\eta})=x^*_{\rm D}(\mathcal{D}_\alpha)$. Therefore, $x^*_{\rm D}(\mathcal{D}_\alpha^{\eta})=x^*_{\rm P}$.

\subsection{Proof of Proposition \ref{pro:omega}}
 Relation (\ref{pro_con_2}) in Theorem \ref{pro:continuous} indicates that $\mathbb{E}[z_{\eta_{(1)}}] \ge \mathbb{E}[z_{\eta_{(2)}}]$ if $\eta_{(1)} < \eta_{(2)}$, which implies that $\mathbb{P}[z_{\eta_{(1)}} > z_{\eta_{(2)}}] \ge \mathbb{P}[z_{\eta_{(1)}} < z_{\eta_{(2)}}]$. Let $(z_{\eta_{(1)}} - z_{\eta_{(2)}})/(\eta_{(1)}-\eta_{(2)})<0$ and $(z_{\eta_{(1)}} - z_{\eta_{(2)}})/(\eta_{(1)}-\eta_{(2)})>0$ be the first and second events respectively (denoted as E1 and E2 respectively), we have $\mathbb{P}$[E1] $\ge$ $\mathbb{P}$[E2]. Given that $\mathcal{D}_\alpha^{\eta_{(1)}}$ and $\mathcal{D}_\alpha^{\eta_{(2)}}$ are two random data sets, it suffices to show that, if $z_{\eta_{(1)}} > z_{\eta_{(2)}}$, $\mathbb{P}[\omega(\eta_{(1)}) > \omega(\eta_{(2)})] \ge \mathbb{P}[\omega(\eta_{(1)}) < \omega(\eta_{(2)})]$.
Further let $(\omega(\eta_{(1)}) - \omega(\eta_{(2)}))/(z_{\eta_{(1)}} > z_{\eta_{(2)}})>0$ and $(\omega(\eta_{(1)}) - \omega(\eta_{(2)}))/(z_{\eta_{(1)}} > z_{\eta_{(2)}}) < 0$ be the third and fourth events respectively (denoted as E3 and E4 respectively), we have $\mathbb{P}$[E3] $\ge$ $\mathbb{P}$[E4]. Finally, considering $(\omega(\eta_{(1)}) - \omega(\eta_{(2)}))/(\eta_{(1)} - \eta_{(2)})<0$ and $(\omega(\eta_{(1)}) - \omega(\eta_{(2)}))/(\eta_{(1)} - \eta_{(2)}) > 0$ as the fifth and sixth events (denoted as E5 and E6 respectively), we have 
 \begin{align}
     \mathbb{P}[\text{E5}]&=\mathbb{P}[\text{E1}]\cdot\mathbb{P}[\text{E3}] + \mathbb{P}[\text{E2}]\cdot\mathbb{P}[\text{E4}] \nonumber \\
      \mathbb{P}[\text{E6}] &= \mathbb{P}[\text{E1}]\cdot\mathbb{P}[\text{E4}] + \mathbb{P}[\text{E2}]\cdot\mathbb{P}[\text{E3}],   \nonumber \\
     \text{and}\;   \mathbb{P}[\text{E5}]-\mathbb{P}[\text{E6}]&=(\mathbb{P}[\text{E1}]-\mathbb{P}[\text{E2}])(\mathbb{P}[\text{E3}]-\mathbb{P}[\text{E4}]) \ge 0. \nonumber
 \end{align}
 Hence, we have we have $\mathbb{P}[\omega(\eta_1) > \omega(\eta_2)]\ge\mathbb{P}[\omega(\eta_1) < \omega(\eta_2)]$ when $\eta_{(1)} < \eta_{(2)}$.


\bibliographystyle{siamplain}
\bibliography{references}
\end{document}